\documentclass[a4paper,pdftex,reqno]{amsart}

\makeatletter
\@namedef{subjclassname@2020}{\textup{2020} Mathematics Subject
Classification}
\makeatother

\usepackage{geometry}
\title{Extensions of the Tong-Yang-Ma representation}
\author{Arthur Souli\'e and Akihiro Takano}


\address{CENTER FOR GEOMETRY AND PHYSICS, INSTITUTE FOR BASIC SCIENCE (IBS), 77, CHEONGAMRO, NAM-GU, POHANG-SI, 790-784, REPUBLIC OF KOREA \& POSTECH, GYEONGSANGBUKDO, REPUBLIC OF KOREA}
\email{arthur.soulie@ibs.re.kr, artsou@hotmail.fr}

\address{GRADUATE SCHOOL OF MATHEMATICAL SCIENCES, THE UNIVERSITY OF TOKYO, 3-8-1 KOMABA, MEGURO-KU, TOKYO, 153-8914, JAPAN}
\email{takano@ms.u-tokyo.ac.jp}

\usepackage{amsmath} 
\usepackage{amsfonts}
\usepackage{mathrsfs}
\usepackage{amssymb}
\usepackage{amsthm}
\usepackage{bm}
\usepackage[new]{old-arrows}
\usepackage{wrapfig}
\usepackage{graphicx}
\usepackage{color}
\usepackage[labelsep=colon]{caption}
\usepackage[labelformat=empty,subrefformat=parens]{subcaption}
\usepackage{multicol}
\usepackage{units}
\usepackage[all]{xy}
\usepackage{tikz}
\usepackage{setspace}
\usepackage{ascmac}
\usepackage{hyperref}
\usepackage[T1]{fontenc}

\newcommand{\relmiddle}[1]{\mathrel{}\middle#1\mathrel{}}

\newtheorem{thm}{Theorem}[section]
\newtheorem{prop}[thm]{Proposition}
\newtheorem{lemma}[thm]{Lemma}
\newtheorem{cor}[thm]{Corollary}

\theoremstyle{definition}
\newtheorem{defi}[thm]{Definition}
\newtheorem{ex}[thm]{Example}
\newtheorem{remark}[thm]{Remark}
\newtheorem{q}[thm]{Question}
\newtheorem{conj}[thm]{Conjecture}

\newcommand{\Z}{\mathbb{Z}}
\newcommand{\R}{\mathbb{R}}
\newcommand{\C}{\mathbb{C}}
\newcommand{\B}{\mathbf{B}}
\newcommand{\WB}{\mathbf{WB}}
\newcommand{\F}{\mathbf{F}}
\newcommand{\Diag}{\mathrm{Diag}}

\begin{document}

\begin{abstract}
In 1996, Tong, Yang and Ma defined a family of representations of the braid group which have the same dimensions as the (unreduced) Burau representations but are not equivalent.
The Burau representation was defined homologically and extended to the string links in several ways.
In this paper, using the method of Silver and Williams, we extend the family of the Tong-Yang-Ma representations to the string links and welded string links. Moreover, we show that the kernels of these representations may be described using some linking numbers. Finally, we apply the Long-Moody construction to the Tong-Yang-Ma representations and study its first properties.
\end{abstract}

\maketitle

\section{Introduction}
Tong, Yang, and Ma \cite{Tong-Yang-Ma} research the representations of the braid group on $n$ strands $\B_n$ such that the $i$-th generator in the Artin presentation maps to the regular matrix $I_{i-1} \oplus T \oplus I_{n-i-1}$, where $I_{k}$ is the $k \times k$ identity matrix and $T$ is an $m \times m$ regular matrix which entries are elements of $\Z [t^{\pm1}]$.
They prove that there exist three kinds of irreducible representations: the trivial one, the (unreduced) Burau one, and a new $n$-dimensional one.
More precisely, when $m=2$, there essentially exist only two non-trivial representations: one is the unreduced Burau representation, and the other is an irreducible representation, called the Tong-Yang-Ma representation.

Furthermore, considering the tensor product $\C \otimes \Z [t^{\pm 1}] \cong \C [t^{\pm 1}]$ and specializing $t$ to complex values in the Burau and Tong-Yang-Ma representations, we obtain complex representations of $\B_n$.
The classification of the irreducible representations of thus type is well researched: by the work of Formanek \cite{Formanek}, for $n \geq 7$, the irreducible complex representations of $\B_n$ of degree $\leq n-1$ are either one-dimensional representation or a tensor product of a one-dimensional representation and a composition factor of the specialization of the (reduced) Burau representation.
Also, we know from Sysoeva \cite{Sysoeva1} that, for $n \geq 9$, the irreducible complex representations of $\B_n$ of degree $n$ are equivalent to a tensor product of one-dimensional representation and specialization of the Tong-Yang-Ma representation.
When $5\leq n\leq 8$, they are classified by Formanek, Lee, Sysoeva and Vazirani \cite{Formanek-Lee-Sysoeva-Vazirani}.
In the case of degree $n+1$, it is shown that there are no irreducible complex representations for $n \geq 10$ by Sysoeva \cite{Sysoeva2}.
Hence, the Burau and Tong-Yang-Ma representations are the only two families of the irreducible complex representations for the braid groups for large $n$ and dimension less or equal to $n+1$.

The problem whether the Burau representation is faithful or not has been highly studied.
For $n = 3$, it is well known to be faithful; see \cite[Theorem 3.15]{Birman}.
On the other hand, Moody \cite{Moody} proved that it is not faithful for $n \geq 9$, and then Long and Paton \cite{Long-Paton} extended that result for $n \geq 6$.
Moreover, Bigelow \cite{Bigelow} showed the non-faithfulness for $n = 5$.
However, for $n = 4$, this problem is still open.
Furthermore, Church and Farb \cite{Church-Farb} show that the kernel of the Burau representation is not finitely generated for $n \geq 6$.
There are also extensions of the Burau representation to string links.
Le Dimet \cite{Le} defined it by using the Fox calculus, and Lin, while Tian and Wang \cite{Lin-Tian-Wang} recovered that extension with a combinatorial and probabilistic approach.
Also, Kirk, Livingston and Wang \cite{Kirk-Livingston-Wang} homologically define the Gassner representation of the string links, which is ``multi-variable version'' of the Burau representation.
In parallel, Silver and Williams \cite{Silver-Williams} defined the two variable Burau matrix $\mathcal{B} \in \mathrm{GL}_n (\mathbb{Z} [u^{\pm1}, v^{\pm 1}])$ for string links.
The specialization $u = 1, v = t^{-1}$ of $\mathcal{B}$ recovers the definition of Lin, Tian and Wang \cite{Lin-Tian-Wang}.

On the other hand, the Tong-Yang-Ma representation is known to not be faithful for $n \geq 3$ by Blanchet and Marin \cite{Blanchet-Marin} by showing that the commutator subgroup of the pure braid group $\mathbf{P}_n$ is in the kernel.
Also, Massuyeau, Oancea and Salamon \cite{Massuyeau-Oancea-Salamon} proved that its kernel may be described via linking numbers.
Apart from these properties, very few researches have been carried out on the Tong-Yang-Ma representations, especially compared to the Burau representations. 

The aim of this paper is to study these representations in greater details.
In Section~\ref{preliminaries}, we recall the notion of string links and the definitions of the Burau and Tong-Yang-Ma representations.
In Section~\ref{ex}, we extend the Tong-Yang-Ma representations to string links in several ways (see Definitions~\ref{def:two variable Tong-Yang-Ma matrix} and \ref{multi-variable Tong-Yang-Ma matrix}), and determine their kernels (see Theorems~\ref{thm:Ker_1} and \ref{thm:Ker_2}).
The ideas of extensions of this representation are based on \cite{Silver-Williams} and \cite{Massuyeau-Oancea-Salamon}.
Moreover, this representation has an extension to the welded braid group, and thus in Section~\ref{welded}, we take this opportunity to define the Tong-Yang-Ma representation for welded string links (see Definition~\ref{def:welded multi-variable Tong-Yang-Ma matrix}).
We also fully determine its kernel using virtual crossings (Theorems~\ref{thm:Ker_3} and \ref{thm:Ker_4}) and reinterpret the welded Tong-Yang-Ma representation in terms of welded biquandles (see Definition~\ref{def:welded Tong-Yang-Ma birack}).
In section~\ref{LM(TYM)}, we apply the Long-Moody construction to the Tong-Yang-Ma representation and prove some first properties.
In particular, it provides a family of irreducible representations for the braid groups (see Theorem~\ref{thm:irred}), which appears to be new as far as we know (see Proposition~\ref{prop:not_equiv_LKB}).
Also the kernels of these representations are strictly included into those of the Burau representations (see Corollary~\ref{coro:kernelLM}).

\section{Preliminaries}\label{preliminaries}

\subsection{String links}
Let $n$ be a positive integer, and $\mathbb{D}^2$ be the unit $2$-disk in $\R^2$.
We fix $n$ distinct points $z_1, \ldots , z_n$ in the interior of $\mathbb{D}^2$.
We may assume that each $z_i$ lies in $\textrm{Int} (\mathbb{D}^2) \cap (\R \times \{ 0 \}) = (-1,1) \times \{ 0 \}$ and that $z_1 < \cdots < z_n$.

\begin{defi}
An \textbf{$n$-string link} is an embedding of the union of $n$ oriented intervals in $\mathbb{D}^2 \times [0,1]$, such that the initial point of each interval corresponds to some $z_i \times \{ 0 \}$ and the endpoint is some $z_j \times \{ 1 \}$.
An interval whose initial point is $z_i \times \{ 0 \}$ is called the \textbf{$i$-th string}.
\end{defi}

\begin{figure}[h]
\begin{center}
\includegraphics[height=130pt]{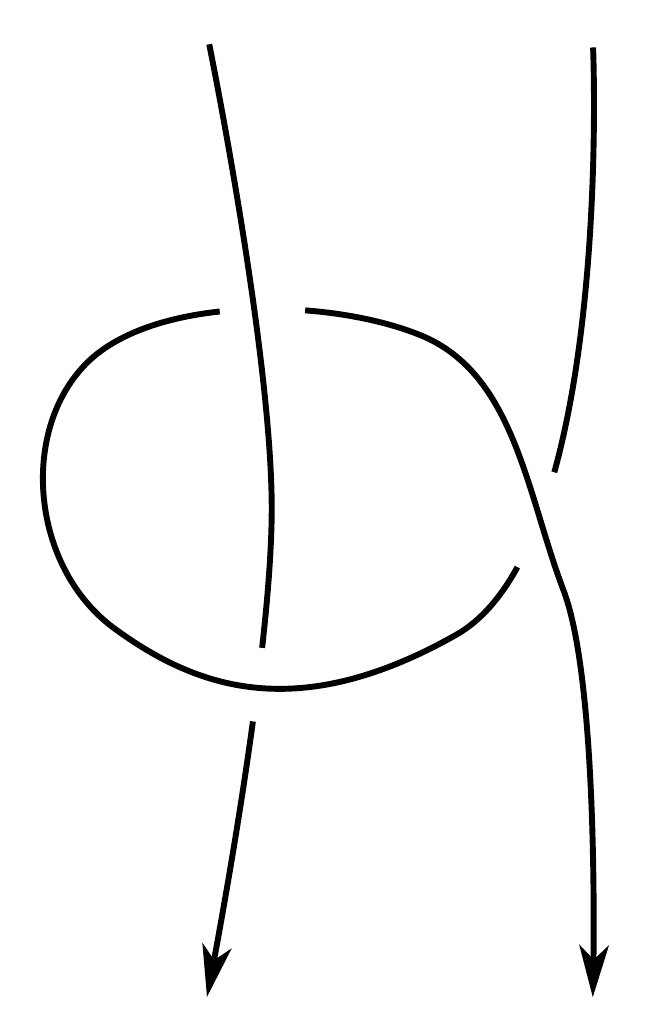}
\caption{A $2$-string link}
\label{string_link1}
\end{center}
\end{figure}

Figure~\ref{string_link1} gives an example of a diagram of a $2$-string link.
A \textbf{pure $n$-string link} is an $n$-string link such that the endpoint of the $i$-th string is $z_i \times \{ 1 \}$ for each $1 \leq i\leq  n$.

\begin{defi}
Two $n$-string links are \textbf{equivalent} if we can pass from one of the associated diagrams to the other one by a finite number of Reidemeister moves (see Figures~\ref{R1} and~\ref{R2R3}) and by an ambient isotopy of the plane.
\end{defi}

\begin{figure}[h]
\begin{center}
\includegraphics[height=105pt]{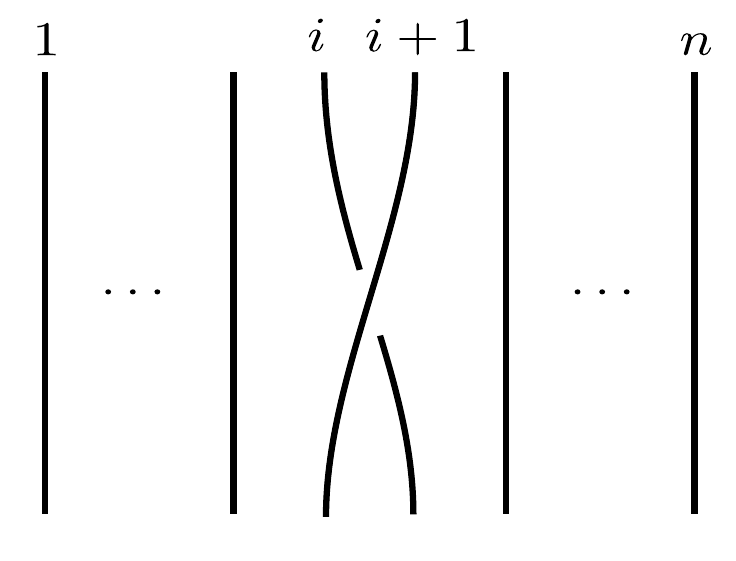}
\caption{Generator $\sigma_i$}
\label{braid_diagram}
\end{center}
\end{figure}

The set of all equivalence classes of $n$-string links, denoted by $\mathcal{SL}_n$, has a structure of a monoid by concatenation and reparametrization.
When each string meets every plane $\mathbb{D}^2 \times \{ t \},\ 0 \leq t \leq 1$, transversely in a single point, the string link is then called an \textbf{$n$-braid}.
The set $\B_n$ of all equivalence classes of $n$-braids is the well-known \textbf{braid group on $n$ strands}. We recall its presentation:
$$
\left \langle \sigma_1 , \ldots , \sigma_{n-1} \relmiddle |
\begin{array}{ll}
\sigma_i \sigma_j = \sigma_j \sigma_i & (|i-j| \geq 2) \\ 
\sigma_i \sigma_{i+1} \sigma_i = \sigma_{i+1} \sigma_i
\sigma_{i+1} & (i = 1 , \ldots , n-2)
\end{array}
\right \rangle.
$$
The $n$-braid in Figure~\ref{braid_diagram} corresponds to the generator $\sigma_i$ of $\B_n$.

\subsection{Tong-Yang-Ma representation} \label{sec:tym}
Tong, Yang, and Ma \cite{Tong-Yang-Ma} classify the representations $\B_n \longrightarrow \mathrm{GL}_{n}(\Z [t^{\pm1}])$ of the form
$$
\sigma_i \longmapsto I_{i-1} \oplus
\left(
\begin{array}{cc}
a & b \\
c & d
\end{array}
\right)
\oplus I_{n-i-1}.
$$
Namely, they prove that there essentially (precisely up to equivalences and transposition) exist only two non-trivial representations of this type.
The first one is the well-known unreduced Burau representation, that is
$$
\sigma_i \longmapsto
I_{i-1} \oplus
\left(
\begin{array}{cc}
0 & t \\
1 & 1-t
\end{array}
\right)
\oplus I_{n-i-1},
$$
that we denote $Bur_n \colon \B_n \longrightarrow \mathrm{GL}_{n} (\Z[t^{\pm1}])$.
The other is the irreducible representation given by
$$
\sigma_i \longmapsto
I_{i-1} \oplus
\left(
\begin{array}{cc}
0 & 1 \\
t & 0
\end{array}
\right)
\oplus I_{n-i-1}.
$$
The later representation is called the \textbf{Tong-Yang-Ma representation}
$$
TYM_n \colon \B_n \longrightarrow \mathrm{GL}_{n} (\Z [t^{\pm1}]).
$$

\section{Extensions of the Tong-Yang-Ma representation} \label{ex}

\subsection{Two variable Tong-Yang-Ma matrix} \label{two}
The procedure of this subsection is based on
\cite{Silver-Williams}.
Let $\Lambda$ be the free abelian group of rank 2 generated by $u$ and $v$.
A \textbf{$\Lambda$-group} is a group $K$ with a right action $K \times \Lambda \rightarrow K$, denoted by $(g,w) \longmapsto g^w$.
Given two $\Lambda$-groups $K_1$ and $K_2$, a map $f \colon K_1 \rightarrow K_2$ is the \textbf{$\Lambda$-homomorphism} if $f$ is a group homomorphism and $f (g^w) = f (g)^w$ for any $g \in K_1$ and $w \in \Lambda$.

\begin{figure}[h]
\begin{center}
\includegraphics[height=140pt]{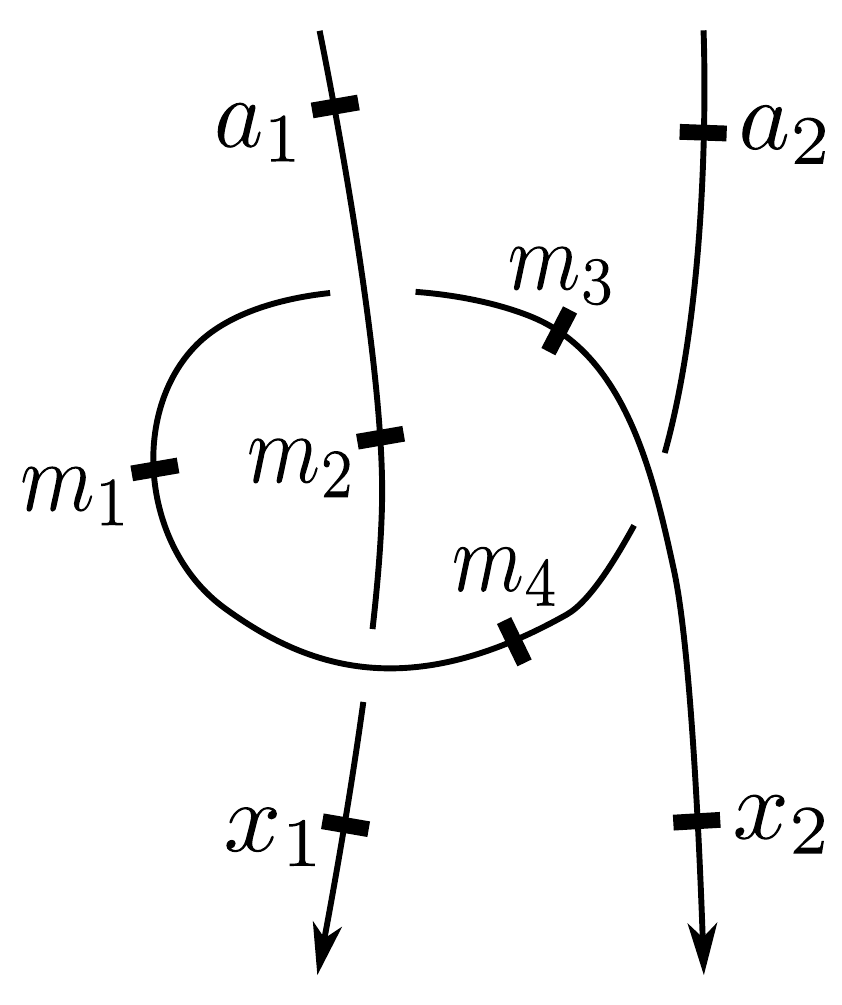}
\caption{A labeled string link}
\label{label_string_link1}
\end{center}
\end{figure}

Let $D$ be a diagram of an $n$-string link $L$.
For each arc of $D$, we put a vertex on it.
Here, we consider each over-crossing arc as the union of two arcs.
If $D$ has $N$ crossings, we need $n+2N$ vertices.
We assume that the top (resp. bottom) vertices are labeled by $a_1, a_2, \ldots, a_n$ (resp. $x_1, x_2, \ldots x_n$) from left to right, and middle vertices, that are neither top nor bottom vertices, are labeled by $m_1, \ldots, m_{2N-n}$ as in Figure~\ref{label_string_link1}.

\begin{figure}[h]
\begin{center}
\includegraphics[height=100pt]{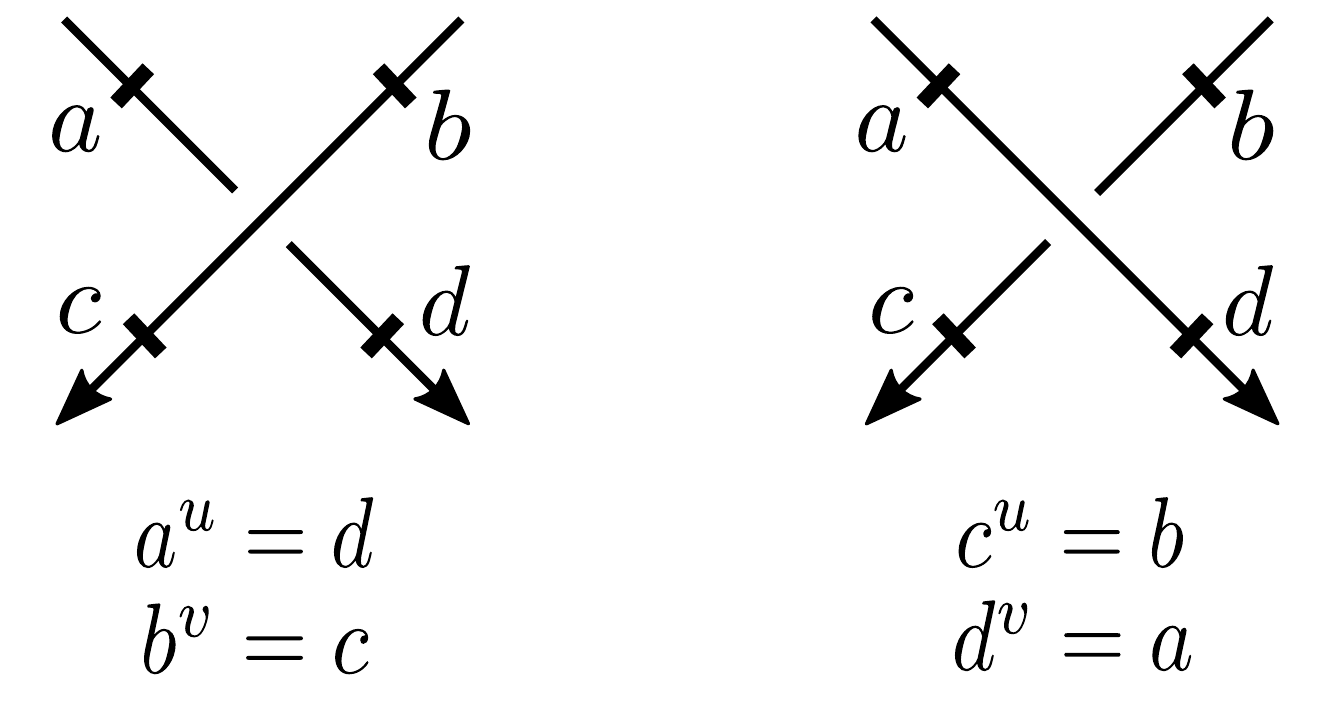}
\caption{The $\Lambda$-relations}
\label{TYM_crossing}
\end{center}
\end{figure}

We construct a $\Lambda$-group associated to a diagram as follows:
\begin{defi}
The $\Lambda$-group $K$ associated to the diagram $D$ is defined by the presentation:
\begin{description}
\item[Generators] the letters indexed by $\mathbb{Z}^{\oplus2}$ corresponding to each label of the vertex of $D$.
For example, for the label $a_1$, the set $\{ (a_1)_{j,k} \}_{j,k}$ is a generating set.
We denote this set by the same letter $a_1$.
The action of $\Lambda$ is defined by $((a_1)_{j,k})^u := a_{j+1,k}$ and $((a_1)_{j,k})^v := a_{j,k+1}$.
\item[Relations] the families of relations indexed by $\mathbb{Z}^{\oplus2}$ corresponding to each crossing of $D$, described in Figure~\ref{TYM_crossing}.
For example, the relations $a^u = d$ and $b^v = c$ denote the families of relations of the form $a_{j+1,k} = d_{j,k}$ and $b_{j,k+1} = c_{j,k}$, respectively, where $j,k \in \mathbb{Z}$.
\end{description}
For a label $a$, the set $a =\{ a_{j,k}, j,k \in \mathbb{Z} \}$ of generators is called the \textbf{$\Lambda$-generator}, and a family of relations is called the \textbf{$\Lambda$-relation}.
\end{defi}

We denote by $\mathcal{S}$ the $n$-tuples of $\Lambda$-generators $ (a_1, \ldots, a_n)$, and by $\mathcal{T}$ the $n$-tuples of $\Lambda$-generators $(x_1, \ldots, x_n)$.

\begin{ex}
The $\Lambda$-group $K$ associated to the diagram in Figure~\ref{label_string_link1} is
$$
K = \left \langle a_1,a_2,m_1,m_2,m_3,m_4,x_1,x_2 \relmiddle | m_1^u=m_3,a_1^v=m_2,m_4^u=a_2,x_2^v=m_3,m_2^u=x_1,m_4^v=m_1 \right \rangle.
$$
Here $\mathcal{S} = (a_1, a_2)$ and $\mathcal{T} = (x_1, x_2)$.
Namely, $K$ has the following presentation:
\begin{align*}
K = \left \langle \left \{
\begin{array}{l}
(a_1)_{j,k}, (a_2)_{j,k}, (m_1)_{j,k}, \\
(m_2)_{j,k}, (m_3)_{j,k}, (m_4)_{j,k}, \\
(x_1)_{j,k}, (x_2)_{j,k}
\end{array}
\right \}_{j,k} \relmiddle | \left \{
\begin{array}{l}
((m_1)_{j,k})^u=(m_3)_{j,k}, ((a_1)_{j,k})^v=(m_2)_{j,k}, \\
((m_4)_{j,k})^u=(a_2)_{j,k}, ((x_2)_{j,k})^v=(m_3)_{j,k}, \\
((m_2)_{j,k})^u=(x_1)_{j,k}, ((m_4)_{j,k})^v=(m_1)_{j,k}
\end{array}
\right \}_{j,k} \right \rangle.
\end{align*}
We can eliminate the generators $m_1,m_2,m_3$ and $m_4$ using the relations, and thus obtain the presentation
$$
K = \left \langle a_1,a_2,x_1,x_2 \relmiddle | x_1 = a_1^{uv},x_2 = a_2 \right \rangle.
$$
In general, for any diagram, the associated $\Lambda$-group $K$ admits a presentation which generators are only the top and bottom vertices.
\end{ex}

Although we may hope that two $\Lambda$-groups associated to equivalent diagrams of the same string link are isomorphic, it is actually not the case. Indeed, a natural question is to check whether or not a $\Lambda$-group is an invariant under the Reidemeister moves R1, R2 and R3 (see Figures~\ref{R1} and \ref{R2R3}) and we prove:

\begin{lemma}
The $\Lambda$-group $K$ associated to the diagram of a string link is invariant under \textup{R2} and \textup{R3}, but is not invariant under \textup{R1}.
\end{lemma}

\begin{proof}
We label the vertices on each diagram of the moves R1, R2 and R3. The $\Lambda$-relations of the right-hand side of R1 (see Figure~\ref{R1}) are the form $b^u = c$  and $a^v = b$.
Hence $c = a^{uv}$ and thus $c \neq a$. Therefore, $K$ is not invariant under R1.

\begin{figure}[h]
\begin{center}
\includegraphics[height=80pt]{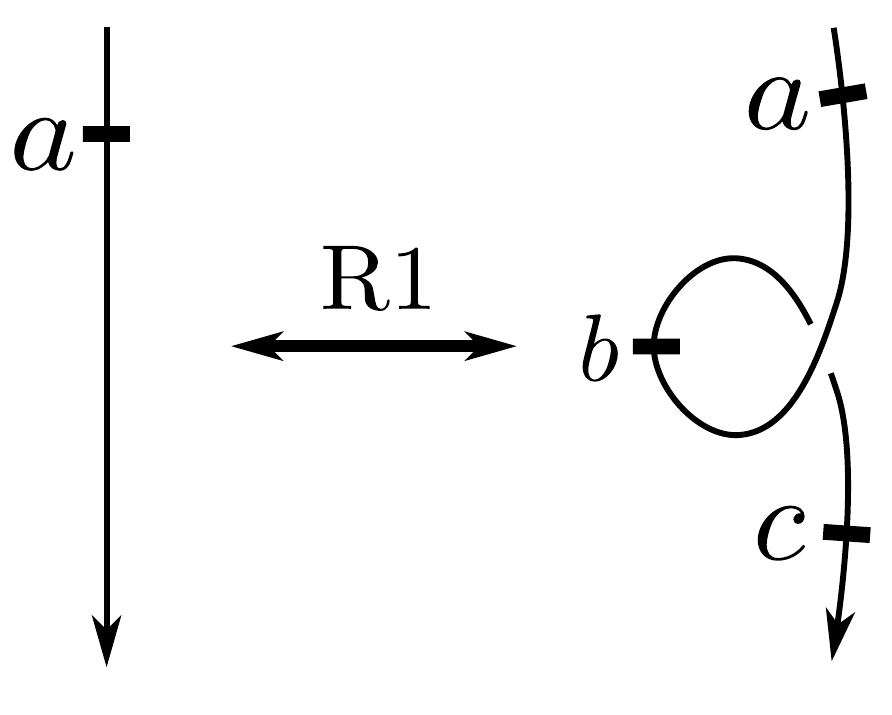}
\caption{The R1-move}
\label{R1}
\end{center}
\end{figure}

\begin{figure}[h]
\begin{tabular}{cc}
\hspace{20pt}
\begin{minipage}{0.33\hsize}
\begin{center}
\includegraphics[height=80pt]{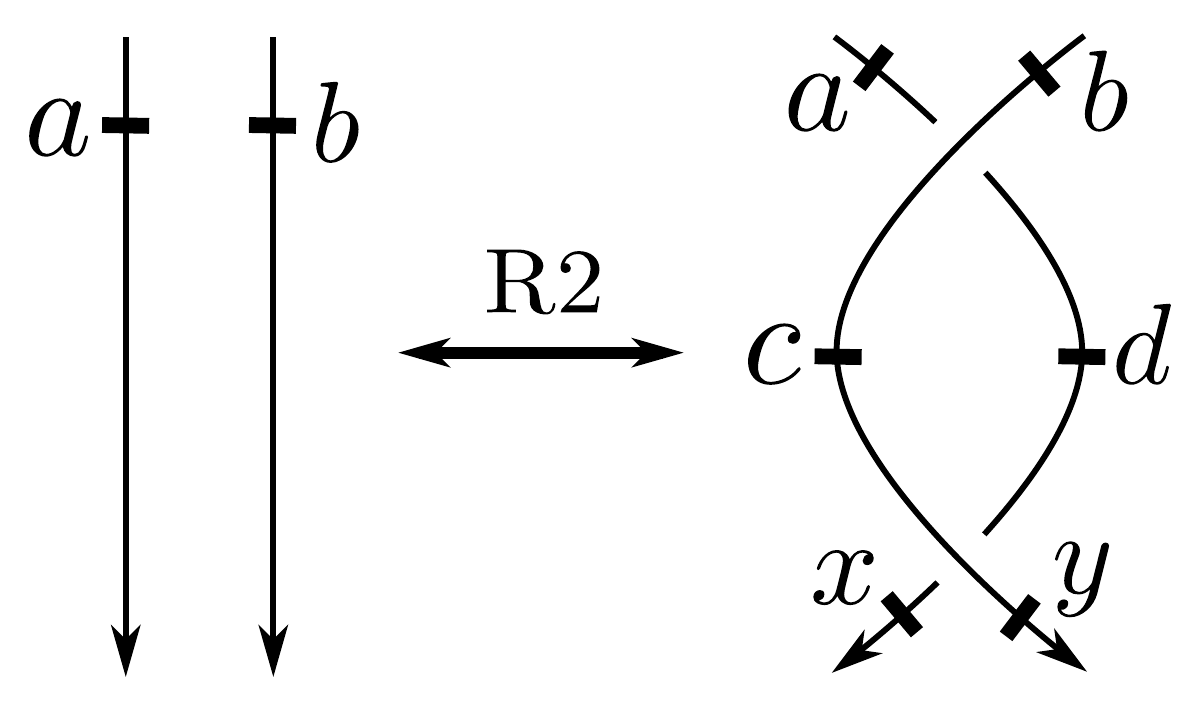}
\subcaption{}
\end{center}
\end{minipage}&
\hspace{30pt}
\begin{minipage}{0.33\hsize}
\begin{center}
\includegraphics[height=80pt]{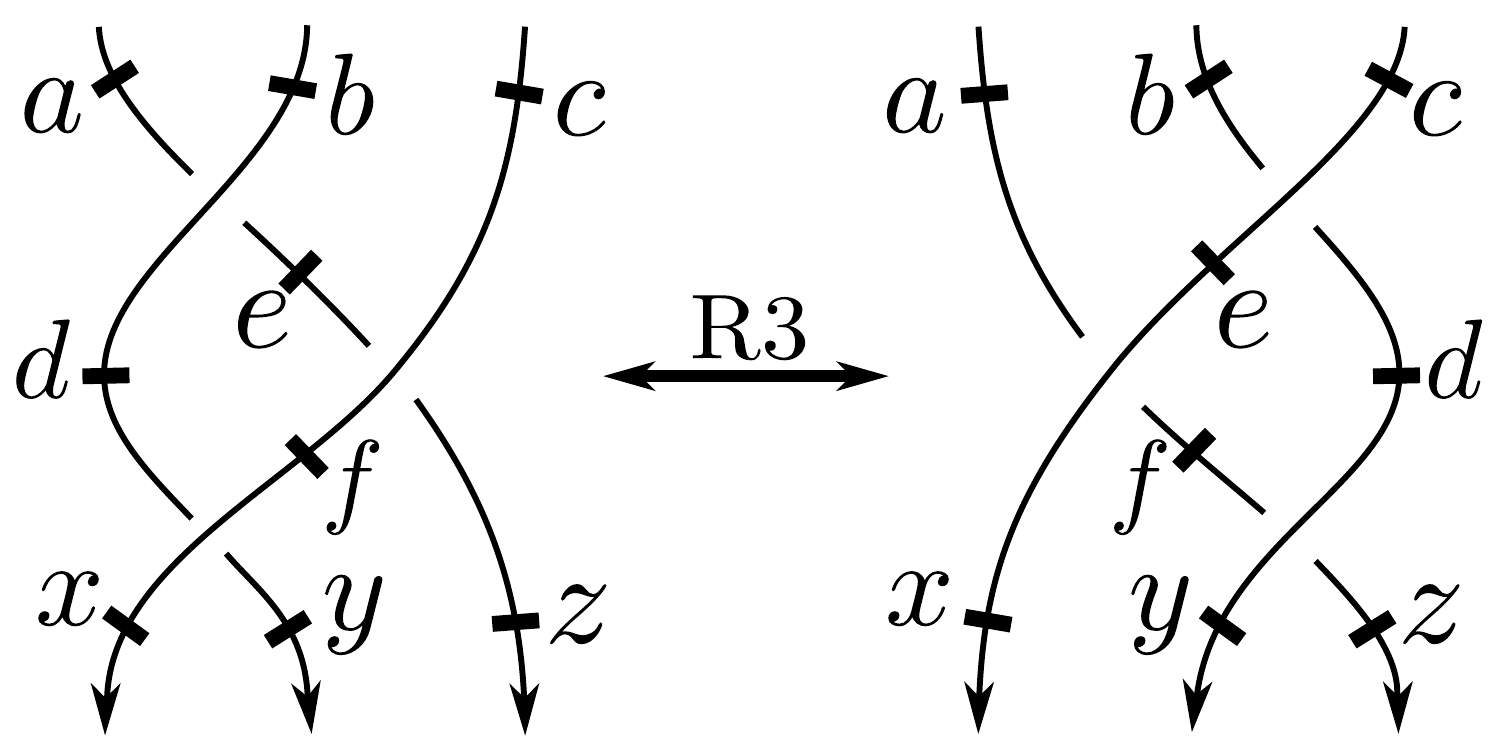}
\subcaption{}
\end{center}
\end{minipage}
\end{tabular}
\vspace{-15pt}
\caption{The R2 and R3 moves}
\label{R2R3}
\end{figure}

The $\Lambda$-relations of the right-hand side of R2 (see the left of Figure~\ref{R2R3}) are the forms
\begin{eqnarray*}
\begin{cases}
a^u = d,\ b^v = c, \\
x^u = d,\ y^v = c.
\end{cases}
\end{eqnarray*}
We deduce that we have the relations $a = x$ and $b = y$. If the orientations of strings of R2 are reversed, we easily check the same relations.

Finally, we prove the invariance of $K$ under R3 (see the right of Figure~\ref{R2R3}).
The invariances under R2 with strings oriented in the same and opposite directions are called the channel unitarity and cross-channel unitarity respectively; see Kauffman \cite{Kauffman}.
In particular, Kauffman \cite{Kauffman} shows that if a quantity has the channel and cross-channel unitarity, then in order to prove the invariance of the quantity under R3 for any orientation of strings it suffices to show the invariance under R3 with strings oriented in the same direction.
The $\Lambda$-relations of the left-hand side and right-hand side of R3 are respectively of the forms
$$
\begin{array}{cc}
\begin{cases}
a^u = e,\ b^v = d, \\
e^u = z,\ c^v = f, \\
d^u = y,\ f^v = x,
\end{cases}
& \hspace{30pt}
\begin{cases}
b^u = d,\ c^v = e, \\
a^u = f,\ e^v = a, \\
f^u = z,\ d^v = y.
\end{cases}
\end{array}
$$
Then, both sets of relations satisfy that $x = c^{v^2},\ y = b^{uv}$ and $z = a^{u^2}$.
\end{proof}

\begin{remark}
In order to extend the Burau representation to string links, Silver and Williams \cite{Silver-Williams} consider the $\Lambda$-relations given by Figure~\ref{Burau_crossing}.
The $\Lambda$-group obtained by using these $\Lambda$-relations is an invariant of the string link.

\begin{figure}[h]
\begin{center}
\includegraphics[height=100pt]{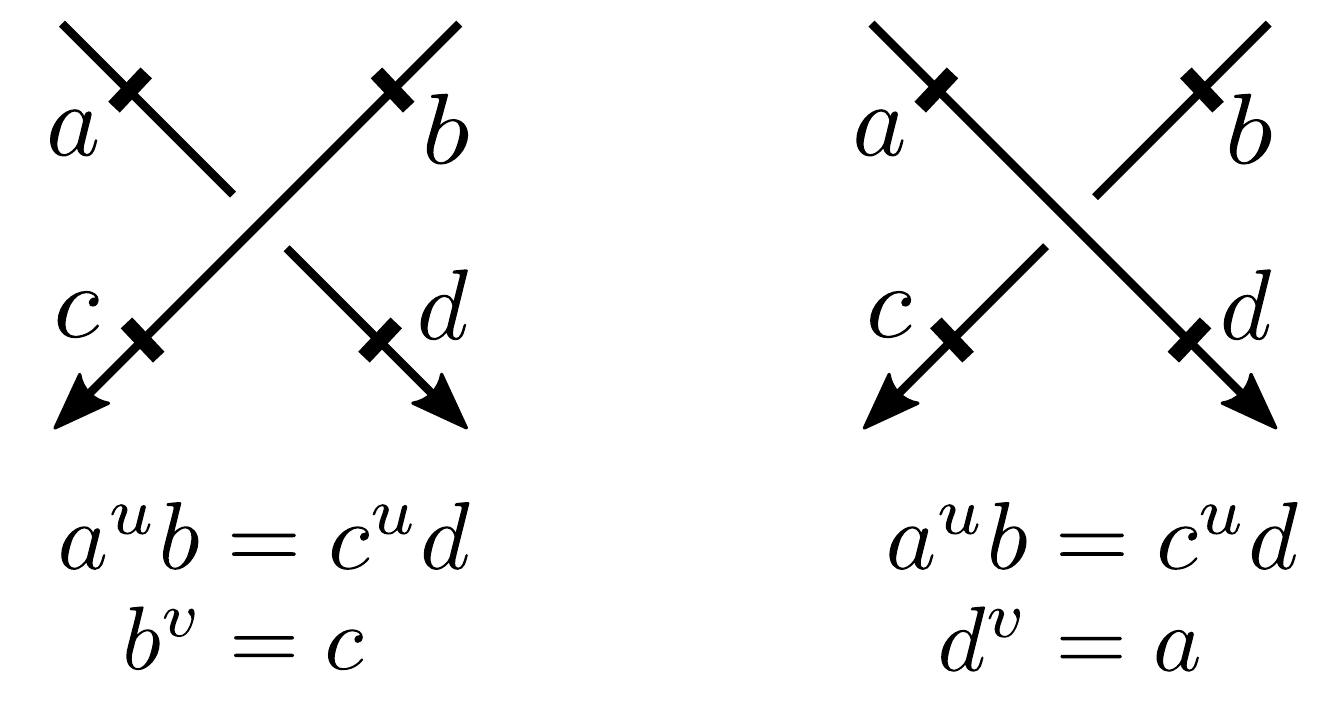}
\caption{$\Lambda$-relations in \cite{Silver-Williams}}
\label{Burau_crossing}
\end{center}
\end{figure}
\end{remark}

We now introduce another $\Lambda$-group associated to $D$ denoted by $\widetilde{K}$, which is also invariant under R1. We may assume without loss of generality that $K$ has a presentation of the form
$$K = \left \langle a_1, \ldots, a_n,x_1, \ldots, x_n \relmiddle | x_1 = a_{i_1}^{w_1}, \ldots, x_n = a_{i_n}^{w_n} \right \rangle,$$
where $w_1, \ldots, w_n \in \Lambda$ correspond to the $\Lambda$-relations of $K$.
Then, we define the $\Lambda$-group $\widetilde{K}$ by
$$
\widetilde{K} := \left \langle a_1, \ldots, a_n,x_1, \ldots, x_n \relmiddle | x_1 = a_{i_1}^{w_1(uv)^{-k_1}}, \ldots, x_n = a_{i_n}^{w_n(uv)^{-k_n}} \right \rangle,
$$
where $k_j$ is the sum of signs of self-crossings in $i_j$-th string. That $\widetilde{K}$ is invariant under three Reidemeister moves being obvious, we deduce that:

\begin{prop}
The $\Lambda$-group $\widetilde{K}$ associated to the diagram of an $n$-string link $L$ is an invariant of $L$.
\end{prop}

The abelianization $\widetilde{K}^{\textrm{ab}}$ of the group $\widetilde{K}$ is a finitely generated $\Lambda$-module.
Moreover, since the associated $\Lambda$-group $\widetilde{K}$ may be with the above defining presentation, the $\Lambda$-module $\widetilde{K}^{\textrm{ab}}$ is free and the sets $\mathcal{S}$ and $\mathcal{T}$ define an ordered basis of $\widetilde{K}^{ab}$.
Hence, we may consider the change of basis matrix between $\mathcal{S}$ and $\mathcal{T}$ over $\mathbb{Z}[\Lambda]$.

\begin{defi}\label{def:two variable Tong-Yang-Ma matrix}
The \textbf{two variable Tong-Yang-Ma matrix} of an $n$-string link $L$ is the change of basis matrix $\mathcal{TYM} (L)$ such that $\mathcal{S}(\mathcal{TYM} (L)) = \mathcal{T}$.
\end{defi}

By construction, we deduce the following properties.

\begin{thm}
\textup{(1)} The matrix $\mathcal{TYM} (L)$ is an invariant of $L$.

\textup{(2)} For any two $n$-string links $L$ and $L'$, we have
$$
\mathcal{TYM} (LL') = \mathcal{TYM} (L) \mathcal{TYM} (L').
$$
In particular, the map
$$
\mathcal{TYM} \colon \mathcal{SL}_n \longrightarrow \mathrm{GL}_n (\mathbb{Z}[\Lambda]) {;}\ L \longmapsto \mathcal{TYM} (L)
$$
is a representation of the $n$-string link monoid $\mathcal{SL}_n$.
\end{thm}
The relations determined by the crossing described in Figure~\ref{TYM_crossing} are of the form
\begin{eqnarray}
\begin{cases}
ua = d,\ vb = c, \hspace{10pt} &\textrm{if} \ \varepsilon = +1\\
u^{-1}b = c,\ v^{-1}a = d, \hspace{10pt} &\textrm{if} \ \varepsilon = -1,
\end{cases}
\end{eqnarray}
where $\varepsilon = \pm 1$ is the sign of the crossing.
The relations (1) can be described in the matrix form $(a,b) M^{\varepsilon} = (c,d)$, where
$$
M =
\left(
\begin{array}{cc}
0 & u \\
v & 0
\end{array}
\right).
$$
The braid group $\B_n$ is generated by $\sigma_1, \ldots, \sigma_{n-1}$.
Then, the matrix $\mathcal{TYM} (\sigma_i)$ has the form:
$$
\mathcal{TYM} (\sigma_i) = I_{i-1} \oplus M \oplus I_{n-i-1}.
$$
When we assign $u = 1$ and $v = t$, we obtain the matrix
$$
I_{i-1} \oplus
\left(
\begin{array}{cc}
0 & 1 \\
t & 0
\end{array}
\right)
\oplus I_{n-i-1}
$$
which coincides with the Tong-Yang-Ma representation.
In this sense, the representation $\mathcal{TYM} \colon \mathcal{SL}_n \rightarrow \mathrm{GL}_n (\mathbb{Z}[\Lambda])$ may be regarded as an extension to string links of the Tong-Yang-Ma representation.

\begin{ex}
We consider the $2$-string link $L$ described in Figure~\ref{label_string_link1}.
The sum of signs of self-crossings in the first string is 0, and that of the second string is $-1$.
Thus, the associated $\Lambda$-group $\widetilde{K}$ is
$$
\widetilde{K} = \left \langle a_1,a_2,x_1,x_2 \relmiddle | x_1 = a_1^{uv},x_2 = a_2^{uv} \right \rangle.
$$
Moreover $\widetilde{K}^{\textrm{ab}}$ has the presentation
$$
\widetilde{K}^{\textrm{ab}} = \left \langle a_1,a_2,x_1,x_2 \relmiddle | x_1 = uva_1,x_2 = uva_2, Comm \right \rangle,
$$
where $Comm$ is the set of all commutators of generators of $\widetilde{K}$.
Hence, the two variable Tong-Yang-Ma matrix is
$$
\mathcal{TYM} (L) = 
\left(
\begin{array}{cc}
uv & 0 \\
0 & uv
\end{array}
\right).
$$
\end{ex}

\begin{ex}
We consider the braid $\sigma_1 \sigma_2^{-1}$; see Figure~\ref{label_braid}. The associated $\Lambda$-group $K$ has the presentation
$$K = \left \langle a_1,a_2,a_3,m_1,x_1,x_2,x_3 \relmiddle | a_1^u = m_1,a_2^v = x_1,x_2^u = a_3,x_3^v = m_1 \right \rangle.$$
We can simplify the relations and thus obtain the following presentation without the middle generator:
$$K = \left \langle a_1,a_2,a_3,x_1,x_2,x_3 \relmiddle | x_1 = a_2^v, x_2 = a_3^{u^{-1}}, x_3 = a_1^{uv^{-1}} \right \rangle.$$
Since there are no self-crossings in $\sigma_1 \sigma_2^{-1}$, $\widetilde{K}$ has the same presentation.
Then, the two variable Tong-Yang-Ma matrix is
$$\mathcal{TYM} (\sigma_1 \sigma_2^{-1}) = 
\left(
\begin{array}{ccc}
0 & 0 & uv^{-1} \\
v & 0 & 0 \\
0 & u^{-1} & 0
\end{array}
\right).$$
Assigning $u = 1$ and $v = t$, we obtain $TYM_3 (\sigma_1 \sigma_2^{-1})$,
$$TYM_3 (\sigma_1 \sigma_2^{-1}) =
\left(
\begin{array}{ccc}
0 & 0 & t^{-1} \\
t & 0 & 0 \\
0 & 1 & 0
\end{array}
\right).$$
\begin{figure}[h]
\begin{center}
\includegraphics[height=130pt]{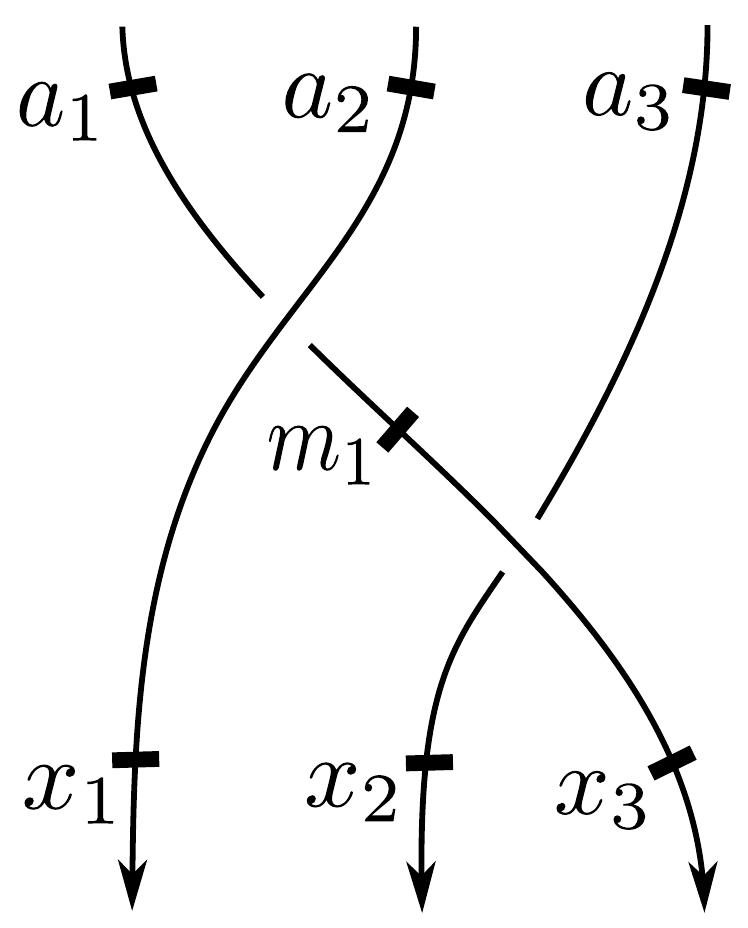}
\caption{Labeled braid $\sigma_1 \sigma_2^{-1}$}
\label{label_braid}
\end{center}
\end{figure}
\end{ex}

\subsection{Multi-variable version}
\begin{figure}[h]
\begin{center}
\includegraphics[height=110pt]{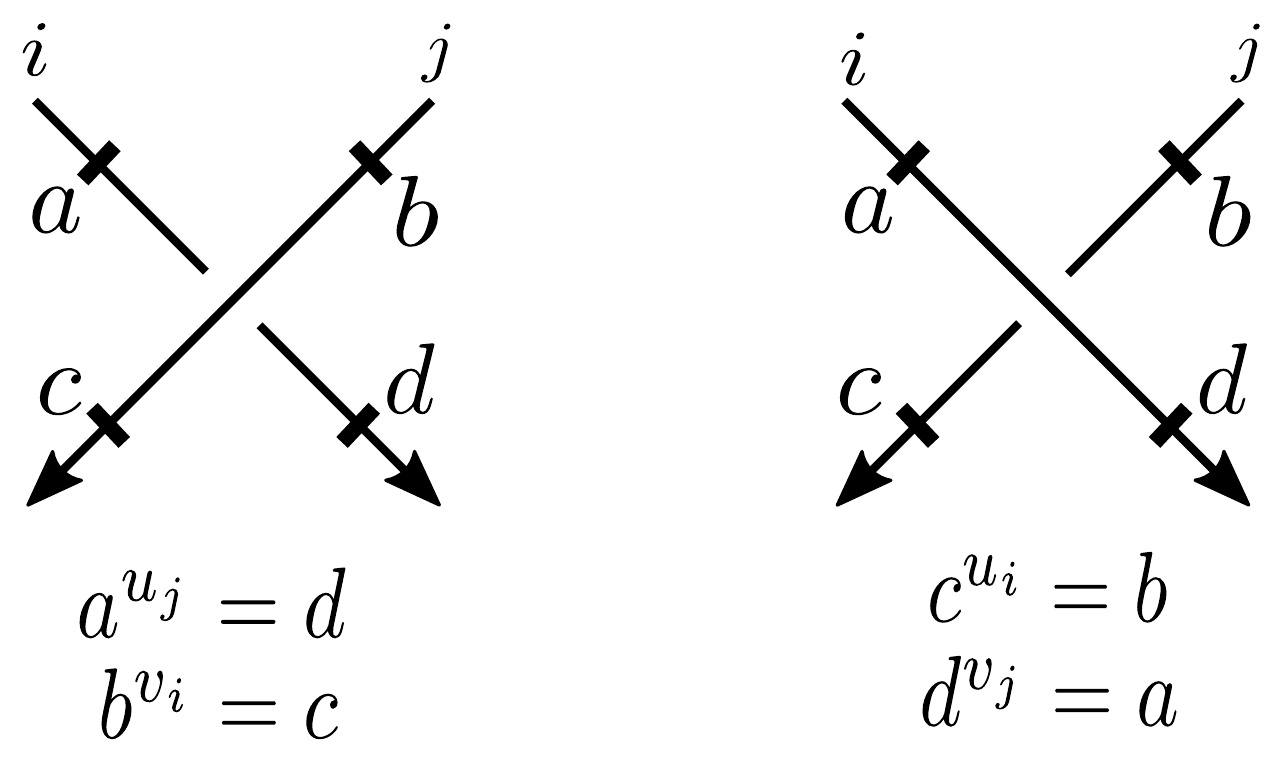}
\caption{$\Lambda_n$-relations}
\label{Multi_TYM_crossing}
\end{center}
\end{figure}

Let $\Lambda_n$ be the free abelian group of rank 2$n$ generated by $u_1, \ldots, u_n, v_1, \ldots, v_n$.
\begin{defi}
Considering a labeled diagram $D$ of a string link $L$, let $K_n$ be the $\Lambda_n$-group defined by the presentation:
\begin{description}
\item[Generators] the letters indexed by $\mathbb{Z}^{\oplus 2n}$ corresponding to each label of the vertex of $D$.
For example, for the label $a_1$, the set $\{ (a_1)_{\bm{n}} \}_{\bm{n} \in \mathbb{Z}^{\oplus 2n}}$ is the generating set.
We denote this set by the same letter $a_1$.
\item[Relations] the families of relations indexed by $\mathbb{Z}^{\oplus2n}$, described in Figure~\ref{Multi_TYM_crossing}.
\end{description}
\end{defi}

As before, we note that the $\Lambda_n$-group $K_n$ is invariant under the Reidemeister moves R2 and R3, but is not invariant under R1. Hence we define the $\Lambda_n$-group $\widetilde{K}_n$ as follows. We may take the top and bottom generators as the generating system of $K_n$, that is,
$$K_n = \left \langle a_1, \ldots, a_n,x_1, \ldots, x_n \relmiddle | x_1 = a_{i_1}^{w_1}, \ldots, x_n = a_{i_n}^{w_n} \right \rangle,$$
where $w_1, \ldots, w_n \in \Lambda_n$ correspond to the $\Lambda_n$-relations of $K_n$.
Then, we define the $\Lambda_n$-group $\widetilde{K}_n$ via the following presentation
$$\widetilde{K}_n := \left \langle a_1, \ldots, a_n,x_1, \ldots, x_n \relmiddle | x_1 = a_{i_1}^{w_1(u_{i_1}v_{i_1})^{-k_1}}, \ldots, x_n = a_{i_n}^{w_n(u_{i_n}v_{i_n})^{-k_n}} \right \rangle.$$
The $\Lambda_n$-group $\widetilde{K}_n$ is an invariant of the $n$-string link $L$, and its abelianization $\widetilde{K}_n^{\textrm{ab}}$ is a free $\Lambda_n$-module of rank $n$ with ordered bases $\mathcal{S}$ and $\mathcal{T}$.

\begin{defi}\label{multi-variable Tong-Yang-Ma matrix}
The \textbf{multi-variable Tong-Yang-Ma matrix} of an $n$-string link $L$ is defined as the change of basis matrix $\mathcal{TYM}_n (L)$ such that $\mathcal{S}(\mathcal{TYM}_n (L)) = \mathcal{T}$.
\end{defi}

We then deduce the following properties on the multi-variable Tong-Yang-Ma matrices:

\begin{thm} \label{mTYM}
\textup{(1)} The matrix $\mathcal{TYM}_n (L)$ is an invariant of $L$.

\textup{(2)} For any $n$-string links $L$ and $L'$, we have
$$
\mathcal{TYM}_n (LL') = \mathcal{TYM}_n (L) (L \cdot \mathcal{TYM}_n (L')).
$$
Here, $\mathcal{SL}_n$ acts on $\Lambda_n$ as follows. If the string of $L$ with endpoint $z_j \times \{ 1 \}$ is the $i_j$-th one, we define the permutation $\tau_L$ by $\tau_L (j) := i_j $ for $1\leq j\leq n$. Then we define the action $\mathcal{SL}_n \times \Lambda_n \rightarrow \Lambda_n$ by $L \cdot u_j := u_{i_j}$ and $L \cdot v_j := v_{i_j}$.
Moreover, this action is extended to $\mathrm{GL}_n (\mathbb{Z}[\Lambda_n])$.
\end{thm}

\begin{proof}
The property \textup{(1)}  is clear from the construction.

For two $n$-string links $L, L'$, we may assume that the associated $\Lambda_n$-groups $\widetilde{K}_n (L), \widetilde{K}_n (L)$ and $\widetilde{K}_n (LL')$ have the following presentations:
\begin{align*}
\widetilde{K}_n (L) &= \left \langle a_1, \ldots, a_n,b_1, \ldots, b_n \relmiddle | b_1 = a_{i_1}^{p_1}, \ldots, b_n = a_{i_n}^{p_n} \right \rangle, \\
\widetilde{K}_n (L') &= \left \langle b_1, \ldots, b_n,c_1, \ldots, c_n \relmiddle | c_1 = b_{j_1}^{q_1}, \ldots, c_n = b_{j_n}^{q_n} \right \rangle, \\
\widetilde{K}_n (LL') &= \left \langle a_1, \ldots, a_n,c_1, \ldots, c_n \relmiddle | c_1 = a_{i_{j_1}}^{r_1}, \ldots, c_n = a_{i_{j_n}}^{r_n} \right \rangle,
\end{align*}
where $p_i, q_i, r_i \in \Lambda_n$ for all $1\leq i\leq n$. The $j$-th string in $L'$ is the $\tau_L(j)$-th string in $LL'$. Hence each number $q_i$ becomes $L \cdot q_i$ in $LL'$ and $r_k = p_{j_k} (L \cdot q_k)$. We therefore obtain the relation $\mathcal{TYM}_n (LL') = \mathcal{TYM}_n (L) (L \cdot \mathcal{TYM}_n (L'))$, thus proving the property \textup{(2)}.
\end{proof}

\begin{ex}
Let $L$ be the $2$-string link described in Figure~\ref{label_string_link1}.
The $\Lambda_2$-group $\widetilde{K}_2$ associated to $L$ has the presentation
$$
\widetilde{K}_2 = \left \langle a_1,a_2,x_1,x_2 \relmiddle | x_1 = a_1^{u_2 v_2}, x_2 = a_2^{u_1 v_1} \right \rangle.
$$
Hence
$$
\mathcal{TYM}_2 (L) =
\left(
\begin{array}{cc}
u_2 v_2 & 0 \\
0 & u_1 v_1
\end{array}
\right).
$$
\end{ex}

\begin{ex}
Consider the braid $\sigma_1 \sigma_2^{-1}$; see Figure~\ref{label_braid}. The associated $\Lambda_3$-group $K_3$ has the presentation
\begin{align*}
K_3 &= \left \langle a_1,a_2,a_3,m_1,x_1,x_2,x_3 \relmiddle | a_1^{u_2} = m_1,a_2^{v_1} = x_1,x_2^{u_1} = a_3,x_3^{v_3} = m_1 \right \rangle \\
&= \left \langle a_1,a_2,a_3,x_1,x_2,x_3 \relmiddle | x_1^u = a_2^{v_1}, x_2 = a_3^{u_1^{-1}}, x_3 = a_1^{u_2v_3^{-1}} \right \rangle \\
&= \widetilde{K}_3.
\end{align*}
Then, the multi-variable Tong-Yang-Ma matrix is
$$
\mathcal{TYM}_3 (\sigma_1 \sigma_2^{-1}) = 
\left(
\begin{array}{ccc}
0 & 0 & u_2v_3^{-1} \\
v_1 & 0 & 0 \\
0 & u_1^{-1} & 0
\end{array}
\right).
$$
Here, $\mathcal{TYM}_3 (\sigma_1)$ and $\mathcal{TYM}_3 (\sigma_2^{-1})$ are of the forms
$$
\begin{array}{cc}
\mathcal{TYM}_3 (\sigma_1) = 
\left(
\begin{array}{ccc}
0 & u_2 & 0 \\
v_1 & 0 & 0 \\
0 & 0 & 1
\end{array}
\right),
& \hspace{20pt}
\mathcal{TYM}_3 (\sigma_2^{-1}) = 
\left(
\begin{array}{ccc}
1 & 0 & 0 \\
0 & 0 & v_3^{-1} \\
0 & u_2^{-1} & 0
\end{array}
\right).
\end{array}
$$
The permutation $\tau_{\sigma_1}$ maps 1 to 2, 2 to 1, and 3 to 3.
Hence
$$
\sigma_1 \cdot \mathcal{TYM}_3 (\sigma_2^{-1}) = 
\left(
\begin{array}{ccc}
1 & 0 & 0 \\
0 & 0 & v_3^{-1} \\
0 & u_1^{-1} & 0
\end{array}
\right)
$$
and $\mathcal{TYM}_3 (\sigma_1 \sigma_2^{-1}) = \mathcal{TYM}_3 (\sigma_1) (\sigma_1 \cdot \mathcal{TYM}_3 (\sigma_2^{-1}))$.
\end{ex}

Let $\mathcal{PSL}_n$ be the pure $n$-string link monoid.
Since the permutation of the pure $n$-string link is trivial, we deduce the following result:

\begin{cor}
For any pure $n$-string links $L$ and $L'$, we have
$$
\mathcal{TYM}_n (LL') = \mathcal{TYM}_n (L) \mathcal{TYM}_n (L').
$$
In particular, the map
$$
\mathcal{TYM}_n \colon \mathcal{PSL}_n \longrightarrow \mathrm{GL}_n (\mathbb{Z}[\Lambda_n]) {;}\ L \longmapsto \mathcal{TYM}_n (L)
$$
is a representation of the pure $n$-string link monoid
$\mathcal{PSL}_n$.
\end{cor}

\subsection{Another extension of Tong-Yang-Ma representation}
The underlying method used in this subsection is inspired from \cite{Massuyeau-Oancea-Salamon}.

Let $w_0$ be a fixed point in $\partial \mathbb{D}^2$, $N (z_i)$ a neighborhood of $z_i$ in $\textrm{Int} (\mathbb{D}^2)$ and $w_i$ a fixed point in $\partial N (z_i)$ for $i=1, \ldots, n$.
We may take the neighborhoods $N(z_i)$ so that they are disjoint to each other.
Set $N (Z) := N(z_1) \sqcup \cdots \sqcup N(z_n)$.
Let $L$ be an $n$-string link and $L_i$ be the $i$-th string of $L$.
For each $1\leq i\leq n$, let $c_i$ be a path from $w_0$ to $w_i$ and set $c_i^j := c_i \times \{ j \} \ (j = 0,1)$ and $c := (c_1, \ldots, c_n)$.

We consider a neighborhood $N(L)$ of $L$ such that $N(L) \cap (\mathbb{D}^2 \times \{ j \}) = N(Z) \times \{ j \}$ for $j =0,1$ and a path $\widetilde{L}_i$ on $\partial N(L)$ from $w_i \times \{ 0 \}$ to $w_{k_i} \times \{ 1 \}$ parallel to $L_i$ such that the linking number of $L_i$ and $\widetilde{L}_i$ is 0, where $k_i$ is the endpoint number of $L_i$.
We set $w := w_0 \times [0,1]$ and orient it from $w_0 \times \{ 0 \}$ to $w_0 \times \{ 1 \}$.
We consider $X := \mathbb{D}^2 \times [0,1] - \textrm{Int}(N(L))$, $X_0 := X \cap (\mathbb{D}^2 \times \{ 0 \})$ and let $\iota \colon X_0 \hookrightarrow X$ be an inclusion map.
Then we have $H_1(X_0; \mathbb{Z}) \cong \mathbb{Z}^{\oplus n} \cong \left \langle t_1, \ldots, t_n \relmiddle | t_i t_j = t_j t_i (\forall i, j) \right \rangle$ and $\iota$ induces an isomorphism $\iota_{*} \colon H_1(X_0; \mathbb{Z}) \longrightarrow H_1(X; \mathbb{Z})$, where $t_i$ is a small loop encircling the point $z_i$ in $X_0$.
Note that since $X_0 = (\mathbb{D}^2 - \textrm{Int}(N(Z))) \times \{ 0 \}$, it does not depend on the string link $L$.

\begin{defi} [cf. {\cite[Section 5]{Massuyeau-Oancea-Salamon}}]
The \textbf{multi-variable Tong-Yang-Ma map} $\mathcal{T}_c \colon \mathcal{SL}_n \rightarrow \mathrm{GL}_n (\mathbb{Z}[H_1(X_0; \mathbb{Z})]) \cong \mathrm{GL}_n (\mathbb{Z}[t^{\pm1}_1, \ldots, t^{\pm1}_n])$ is defined for any $n$-string link $L$ by
$$
(\mathcal{T}_c (L))_{ij} :=
\begin{cases}
\iota^{-1}_{*} (c_i^0 \cdot \widetilde{L}_i \cdot (c_j^1)^{-1} \cdot w^{-1}) & \textrm{if} \ \tau_L (j) = i, \\
0 & \textrm{otherwise.}
\end{cases}
$$
Here, the symbol $\cdot$ denotes the product of paths, and thus $c_i^0 \cdot \widetilde{L}_i \cdot (c_j^1)^{-1} \cdot w^{-1}$ is an element of $ H_1 (X ; \mathbb{Z})$.
\end{defi}

\begin{figure}[h]
\begin{center}
\includegraphics[height=200pt]{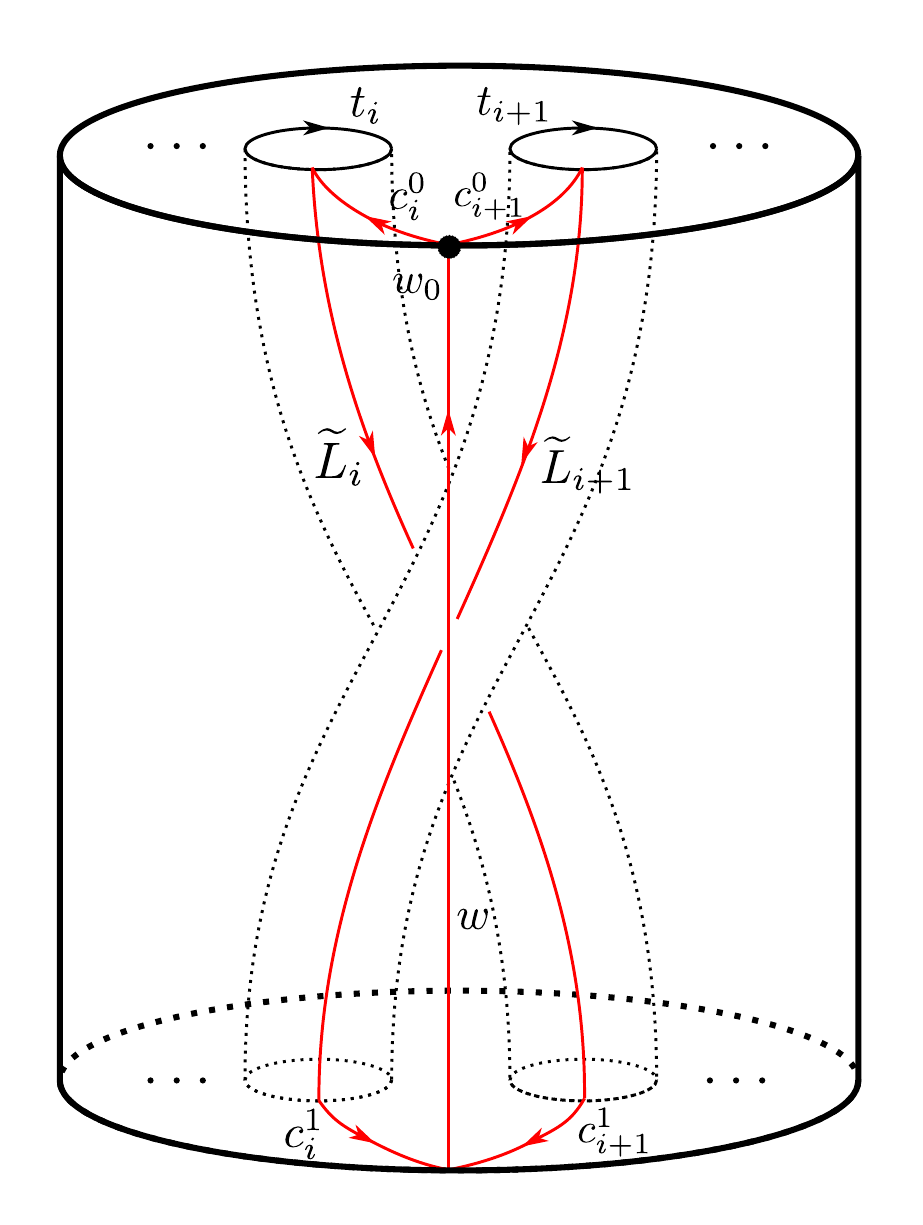}
\caption{Visualization of $\mathcal{T}_c (\sigma_i)$}
\label{braid_homo}
\end{center}
\end{figure}

We depict an example of $\mathcal{T}_c$ in Figure~\ref{braid_homo}.
In this example, since $c_i^0 \cdot \widetilde{L}_i \cdot (c_{i+1}^1)^{-1} \cdot w^{-1} = t_{i+1}$ and $c_{i+1}^0 \cdot \widetilde{L}_{i+1} \cdot (c_i^1)^{-1} \cdot w^{-1} = 1$, we have
$$
\mathcal{T}_c (\sigma_i) =
I_{i-1} \oplus
\left(
\begin{array}{cc}
0 & t_{i+1} \\
1 & 0
\end{array}
\right)
\oplus I_{n-i-1}.
$$

\begin{figure}[h]
\begin{center}
\includegraphics[height=80pt]{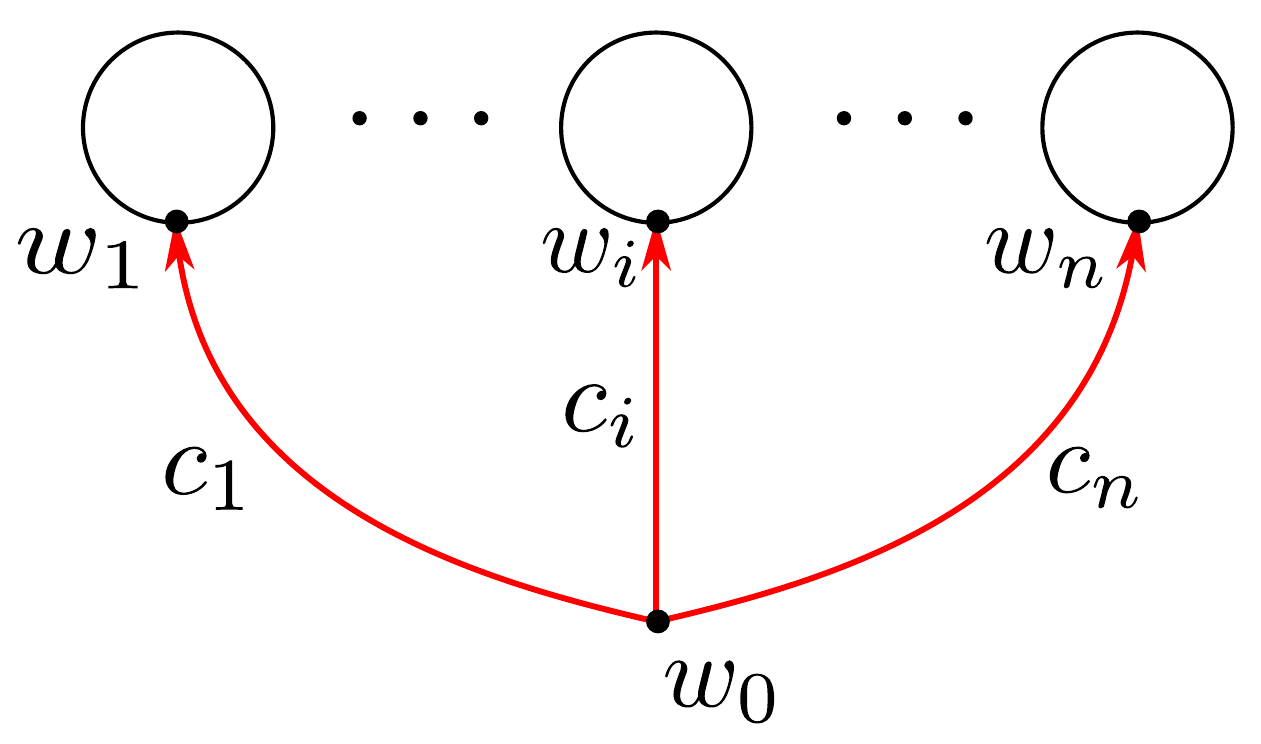}
\caption{\lq\lq Appropriate'' paths}
\label{appropriate_c}
\end{center}
\end{figure}

We consider paths $c = (c_1, \ldots, c_n)$ which are standard in the sense that they are pairwise disjoint and the product $c_{i} \cdot \partial N(z_{i}) \cdot c_{i}^{-1}$ is homotopic to the $i$-th standard generator of the fundamental group of $\mathbb{D}^{2}-\{z_1,...,z_n\}$ ; see Figure~\ref{appropriate_c}. Then $\mathcal{T}_c$ coincides with the specialization $u_1 = t_1, \ldots, u_n = t_n$ and $v_1 = \cdots = v_n = 1$ of the multi-variable Tong-Yang-Ma matrix $\mathcal{TYM}_n$.
Indeed, the diagram of Figure~\ref{braid_homo} corresponds to the left-hand side of Figure~\ref{Multi_TYM_crossing} specialized $u_j = t_j$ and $v_i = 1$.

\begin{ex}
We consider standard paths as above and consider the $2$-string link $L$ in Figure~\ref{string_link_homo} whose diagram coincides with Figure~\ref{string_link1}.
The element in $H_1(X; \mathbb{Z})$ corresponding to the first string is $t_2$, and the one corresponding to the second string is $t_1$. Hence:
$$
\mathcal{T}_c (L) =
\left(
\begin{array}{cc}
t_2 & 0 \\
0 & t_1
\end{array}
\right).
$$
This matrix coincides with the specialization $u_1 = t_1, u_2 = t_2$ and $v_1=v_2=1$ of $\mathcal{TYM}_2 (L)$.
\begin{figure}[h]
\begin{center}
\includegraphics[height=200pt]{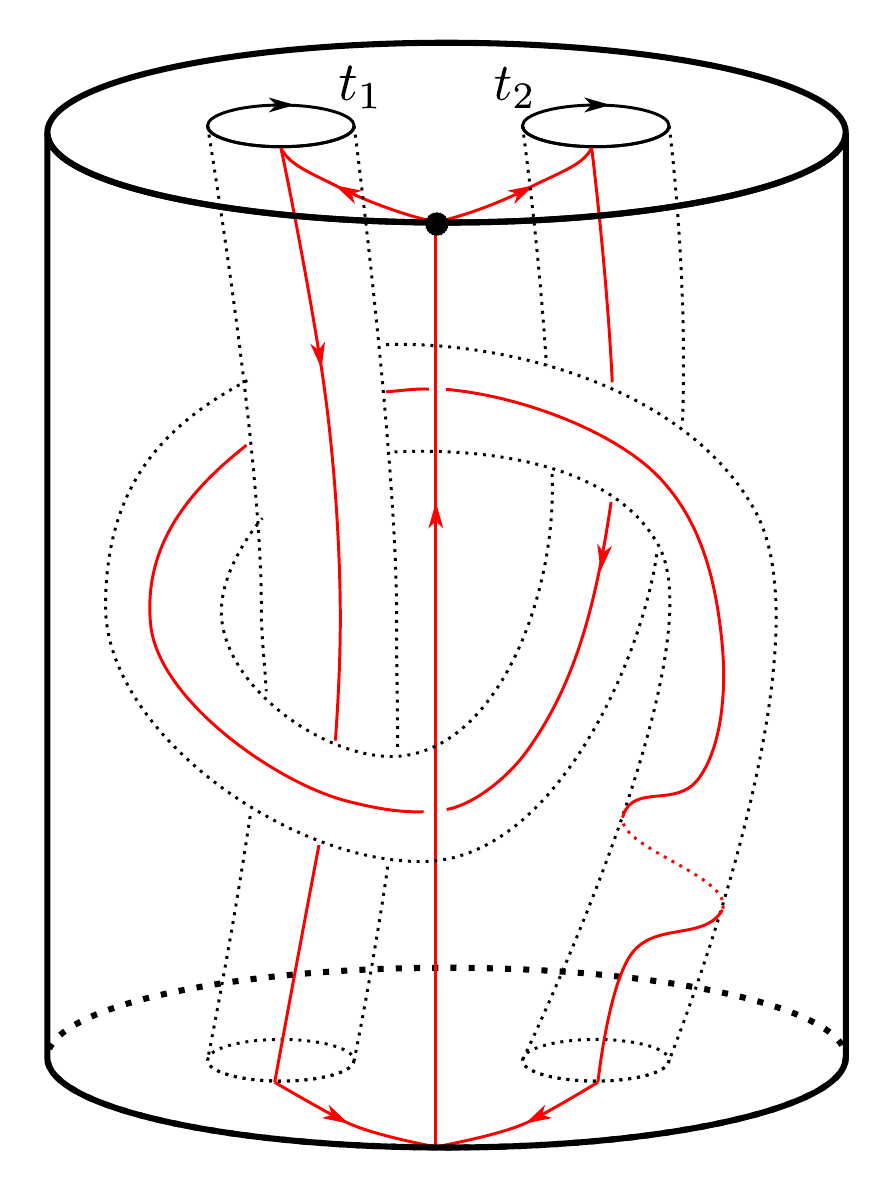}
\caption{$2$-string link $L$}
\label{string_link_homo}
\end{center}
\end{figure}
\end{ex}

\begin{figure}[h]
\begin{center}
\includegraphics[height=110pt]{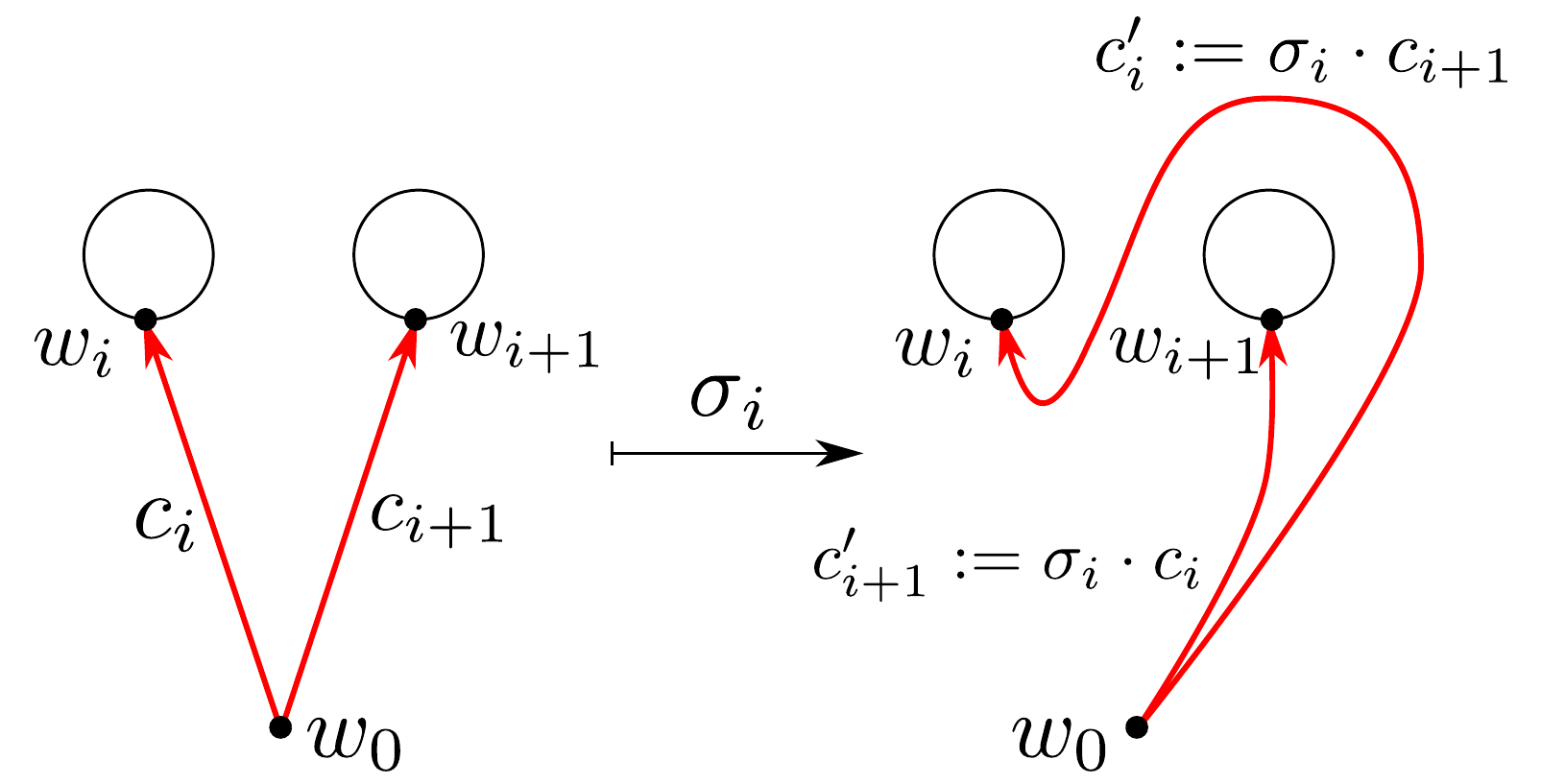}
\caption{Braid action on $c$}
\label{braid_action}
\end{center}
\end{figure}

There is a natural action of the braid group $\B_n$ on these paths $c = (c_1, \ldots, c_n)$, that we denote by $\sigma \cdot c := (\sigma \cdot c_1, \ldots, \sigma \cdot c_n)$.
We set $c' := (c'_1, \ldots, c'_n) := (\sigma \cdot c_{i_1}, \ldots, \sigma \cdot c_{i_n})$, where $i_j$ satisfies $\tau_{\sigma} (i_j) = j \ (j = 1, \ldots, n)$.
Figure~\ref{braid_action} illustrates the action of the generator $\sigma_i$ on $c$.

\begin{thm}[cf. {\cite[Theorem 5.5]{Massuyeau-Oancea-Salamon}}] \label{MOS}

\textup{(1)} For any path $c = (c_1, \ldots, c_n)$ and two $n$-string links $L$ and $L'$, we have
$$
\mathcal{T}_c (LL') = \mathcal{T}_c (L) (L \cdot \mathcal{T}_c (L')),
$$
where the action of $\mathcal{SL}_n$ on $\mathrm{GL}_n (\mathbb{Z}[H_1(X_0; \mathbb{Z})])$ is defined by extending the action of $\mathcal{SL}_n$ on $H_1(X_0; \mathbb{Z})$ given by $L \cdot t_j := t_{\tau_L(j)}$.

\textup{(2)} For any path $c$, any $n$-braid $\sigma$ and any $n$-string link $L$, we have
$$
\mathcal{T}_{\sigma \cdot c} (L) = \sigma_* (\mathcal{T}_c (\sigma)^{-1} \mathcal{T}_c (L) (L \cdot \mathcal{T}_c (\sigma))),
$$
where the action of $\B_n$ on $\mathrm{GL}_n (\mathbb{Z}[t^{\pm1}_1, \ldots, t^{\pm1}_n])$ is defined by $\sigma_* (a_{ij}) := (a_{\tau_{\sigma} (i) \tau_{\sigma} (j)})$ for any $(a_{ij}) \in \mathrm{GL}_n (\mathbb{Z}[t^{\pm1}_1, \ldots, t^{\pm1}_n])$.
\end{thm}

\begin{proof}
(1) We assume that $j = \tau_{L'} (k)$ and $i = \tau_{L} (j) =\tau_{LL'} (k)$, that is, the $i$-th string $L_i$ of $L$ joins the $j$-th string $L'_j$ of $L'$.
Hence $L'_j$ is actually the $i$-th string in $LL'$.
Thus we need to correct which strings are tangled with the string $L'_j$ by using $\tau_{L}$.
We deduce that:
$$
(\mathcal{T}_c (LL'))_{ik} = (\mathcal{T}_c (L))_{ij} (L \cdot (\mathcal{T}_c (L'))_{jk}).
$$
If $i \neq \tau_{LL'} (k)$, then the $(i,k)$-entries of both matrices $\mathcal{T}_c (LL')$ and $\mathcal{T}_c (L) (L \cdot \mathcal{T}_c (L'))$ are zero.
Therefore, we obtain the equality.

(2) We consider the $n$-string link $\sigma^{-1} L \sigma$ and assume that $\sigma^{-1}$ is in $\mathbb{D}^2 \times [0,1/3]$, that $L$ is in $\mathbb{D}^2 \times [1/3,2/3]$ and that $\sigma$ is in $\mathbb{D}^2 \times [2/3,1]$.
Set $\tau_L (k) := j$.
Then
\begin{align*}
(\mathcal{T}_{\sigma \cdot c} (L))_{jk} &= (\sigma \cdot c_{i_j}^0) \cdot \widetilde{L}_j \cdot (\sigma \cdot c_{i_k}^1)^{-1} \cdot w^{-1} \\
&= (c_{i_j}^0 \cdot \widetilde{\sigma^{-1}} \cdot (c_j^{1/3})^{-1}) \cdot (c_j^{1/3} \cdot \widetilde{L}_j \cdot (c_k^{2/3})^{-1}) \cdot (c_k^{2/3} \cdot \widetilde{\sigma} \cdot (c_{i_k}^1)^{-1}) \cdot w^{-1},
\end{align*}
where $\tau_{\sigma}(i_j) = j$ and $\tau_{\sigma}(i_k) = k$.
Here, we compute that
\begin{align*}
& c_{i_j}^0 \cdot \widetilde{\sigma^{-1}} \cdot (c_j^{1/3})^{-1} \cdot w|_{[0,1/3]}^{-1} = \sigma \cdot (\mathcal{T}_c (\sigma^{-1}))_{i_{j}j}, \\
& c_j^{1/3} \cdot \widetilde{L}_j \cdot (c_k^{2/3})^{-1} \cdot w|_{[1/3,2/3]}^{-1} = (\mathcal{T}_c (L))_{jk}, \\
& c_k^{2/3} \cdot \widetilde{\sigma} \cdot (c_{i_k}^1)^{-1} \cdot w|_{[2/3,1]}^{-1} = L \cdot (\mathcal{T}_c (\sigma))_{ki_k}.
\end{align*}
Hence
$$
(\mathcal{T}_{\sigma \cdot c} (L))_{jk} = (\sigma \cdot (\mathcal{T}_c (\sigma^{-1}))_{i_jj}) (\mathcal{T}_c (L))_{jk} (L \cdot (\mathcal{T}_c (\sigma))_{ki_k}).
$$
Therefore, we have
$$
\mathcal{T}_{\sigma \cdot c} (L) = \sigma_* ((\sigma \cdot \mathcal{T}_c (\sigma^{-1})) \mathcal{T}_c (L) (L \cdot \mathcal{T}_c (\sigma))).
$$
Moreover, since
$I_n = \mathcal{T}_c (\sigma \sigma^{-1}) = \mathcal{T}_c (\sigma) (\sigma \cdot \mathcal{T}_c (\sigma^{-1}))$,
we deduce that $\sigma \cdot \mathcal{T}_c (\sigma^{-1}) = \mathcal{T}_c (\sigma)^{-1}$, which ends the proof.
\end{proof}

\begin{cor}
\textup{(1)} For any path $c$ and two pure $n$-string links $L, L'$, we have
$$
\mathcal{T}_c (LL') = \mathcal{T}_c (L) \mathcal{T}_c (L').
$$

\textup{(2)} For any path $c$, any $n$-braid $\sigma$ and any pure $n$-string link $L$, we have
$$
\mathcal{T}_{\sigma \cdot c} (L) = \mathcal{T}_c (L).
$$
\end{cor}

\begin{proof}
The property (1) is clear. For the property (2), we deduce from Theorem~\ref{MOS} (2) that:
$$
\mathcal{T}_{\sigma \cdot c} (L) = \sigma_* (\mathcal{T}_c (\sigma)^{-1} \mathcal{T}_c (L) \mathcal{T}_c (\sigma)).
$$
Here, $\mathcal{T}_c (L)$ is a diagonal matrix, and so is $\mathcal{T}_c (\sigma)^{-1} \mathcal{T}_c (L) \mathcal{T}_c (\sigma)$.
Let $a_i$ be the $i$-th diagonal entry of $\mathcal{T}_c (L)$, then the $j$-th diagonal entry of $\mathcal{T}_c (\sigma)^{-1} \mathcal{T}_c (L) \mathcal{T}_c (\sigma)$ is $a_{i_j}$, where $\tau_{\sigma} (i_j) = j$ for $j = 1, \ldots, n$.
Hence:
$$
\mathcal{T}_c (\sigma)^{-1} \mathcal{T}_c (L) \mathcal{T}_c (\sigma) = \textrm{Diag}(a_{\tau_{\sigma}^{-1} (1)}, \ldots, a_{\tau_{\sigma}^{-1} (n)}) = \sigma_*^{-1} \mathcal{T}_c (L).
$$
Therefore we have $\mathcal{T}_{\sigma \cdot c} (L) = \sigma_*( \sigma_*^{-1} \mathcal{T}_c (L)) = \mathcal{T}_c (L)$.
\end{proof}

\subsection{Kernel of Tong-Yang-Ma representation} \label{ker}
We consider the two variable Tong-Yang-Ma representation $\mathcal{TYM} \colon \mathcal{SL}_n \rightarrow \mathrm{GL}_n (\mathbb{Z}[\Lambda])$.
It clearly follows from the definitions that the $n$-string link $L$ is pure if and only if $\mathcal{TYM} (L)$ is the diagonal matrix.
Moreover, assuming that the $\Lambda$-group $\widetilde{K}$ associated to $L$ has the presentation
$$\widetilde{K} = \left \langle a_1, \ldots, a_n,x_1, \ldots, x_n \relmiddle | x_1 = a_1^{w_1}, \ldots, x_n = a_n^{w_n} \right \rangle,$$
where $w_1, \ldots, w_n \in \Lambda$.
Then $\mathcal{TYM} (L) = \textrm{Diag} (w_1, \ldots, w_n)$.
For each $1\leq i \leq n$, we may set $w_i = (uv)^{m_i} \ (m_i \in \mathbb{Z})$ since $L$ is pure.
By $\Lambda$-relations, $m_i$ is equal to the sum of signs of under (over)-crossings between the $i$-th string and the other string.
Since we canceled the sum of signs of self-crossings by $-k_i$, 2$m_i$ is the sign of all crossings of $i$-th and the other strings, and thus
$$
m_i = \sum_{j = 1, j \neq i}^{n} \ell_{ij} (L),
$$
where $\ell_{ij} (L)$ is the linking number of the $i$-th and $j$-th strings of $L$.
Therefore, we obtain the following result.

\begin{thm}\label{thm:Ker_1}
For $n \geq 2$, we have
$$
\ker \mathcal{TYM} =
\left \{ L \in \mathcal{PSL}_n \relmiddle | 
{\rm for \ each}\ 1\leq i \leq n, \ \sum_{j = 1, j \neq
i}^{n} \ell_{ij} (L) = 0 \right \}.
$$
\end{thm}

When we consider the multi-variable representation $\mathcal{TYM}_n \colon \mathcal{PSL}_n \to \mathrm{GL}_n (\mathbb{Z}[\Lambda_n])$, the $\Lambda_n$-relations distinguish which strings get tangled.
Hence its kernel may be written as follows.

\begin{thm}\label{thm:Ker_2}
For $n \geq 2$,
$$
\ker \mathcal{TYM}_n =
\left \{ L \in \mathcal{PSL}_n \relmiddle | 
{\rm for \ any} \ 1 \leq i \neq j \leq n,\ \ell_{ij} (L) = 0 \right \}.
$$
\end{thm}

\begin{remark}
Massuyeau, Oancea and Salamon \cite{Massuyeau-Oancea-Salamon} introduce the Picard-Lefschetz monodromy cocycle $\mathcal{S} \colon \widetilde{\B}_n \rightarrow \mathrm{GL}_n (\mathbb{Z} [\F_n])$, where $\widetilde{\B}_n$ is the framed braid group.
This cocycle is not a representation.
However, if we restrict it to the braid group and take the composite with the projection $\F_n \rightarrow \mathbb{Z} \cong \langle t \rangle$, then this becomes a representation and coincides with the Tong-Yang-Ma representation.
Also, restricting this cocycle to the pure braid group and composing with the natural projection $\F_n \rightarrow \mathbb{Z}^n \cong \left \langle t_1, \ldots, t_n \relmiddle | t_i t_j = t_j t_i \ (\forall i, j) \right \rangle$, they obtain the multi-variable Tong-Yang-Ma representation $\hat{\mathcal{S}} \colon \mathbf{P}_n \rightarrow \mathrm{GL}_n (\mathbb{Z} [t^{\pm1}_1, \ldots , t^{\pm1}_n])$, where $\mathbf{P}_n$ is the pure braid group.
\cite[Proposition 4]{Massuyeau-Oancea-Salamon} says that
$$
\hat{\mathcal{S}} (b) = \textrm{Diag} \left( \prod_{j \neq 1} t^{-\ell_{1j} (b)}_j , \ldots , \prod_{j \neq n} t^{-\ell_{nj} (b)}_j \right)
$$
for any $b \in \mathbf{P}_n$.
Therefore, specializing $t_1 = \cdots = t_n = t$, the exponent of $t$ of each diagonal entry is the sum of linking numbers.
\end{remark}

\begin{remark}
As above, the Tong-Yang-Ma representation thus has a relationship with the linking number. Hence, if the construction using Silver-Williams's method may seem verbose, this extension can also be used for the virtual string link and the welded string link; see Section~\ref{welded}.
\end{remark}

\section{Welded version}\label{welded}

\subsection{Definitions}
We fix $n$ real numbers $0 < z_1 < \cdots < z_n < 1$.

\begin{figure}[h]
\begin{tabular}{ccc}
\begin{minipage}{0.33\hsize}
\begin{center}
\includegraphics[height=50pt]{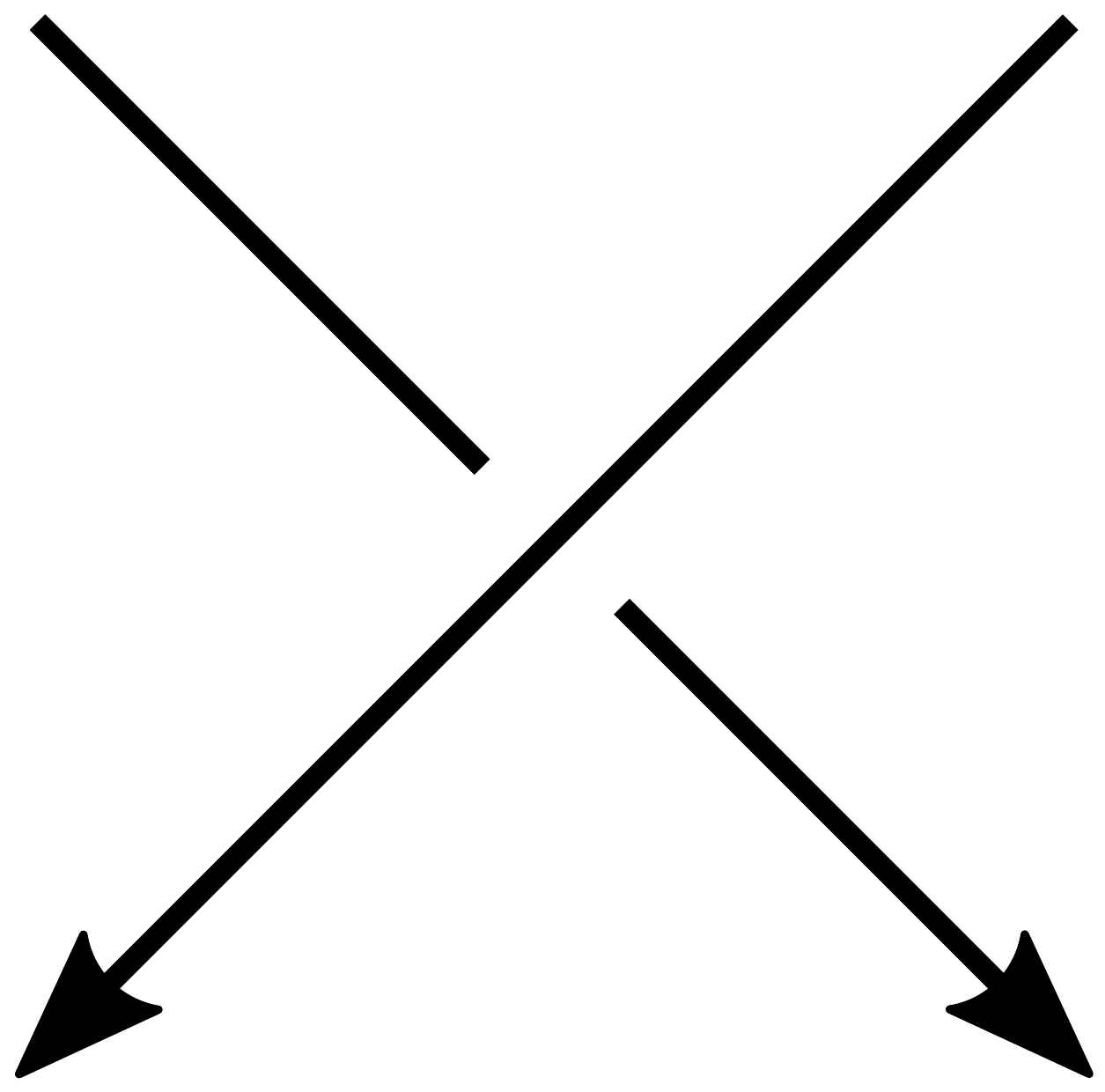}
\subcaption{Positive crossing}
\end{center}
\end{minipage}
\begin{minipage}{0.33\hsize}
\begin{center}
\includegraphics[height=50pt]{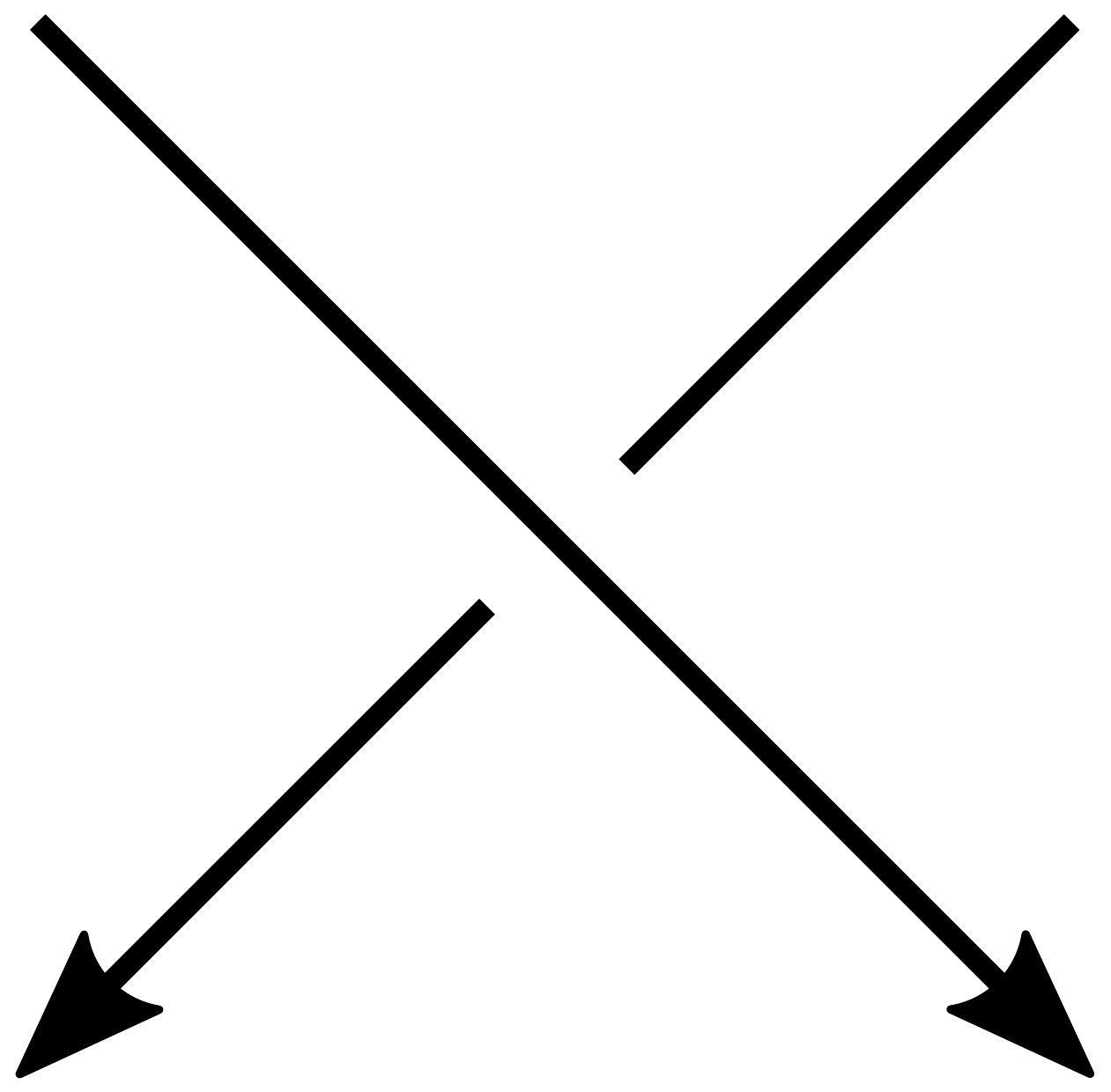}
\subcaption{Negative crossing}
\end{center}
\end{minipage}
\begin{minipage}{0.33\hsize}
\begin{center}
\includegraphics[height=50pt]{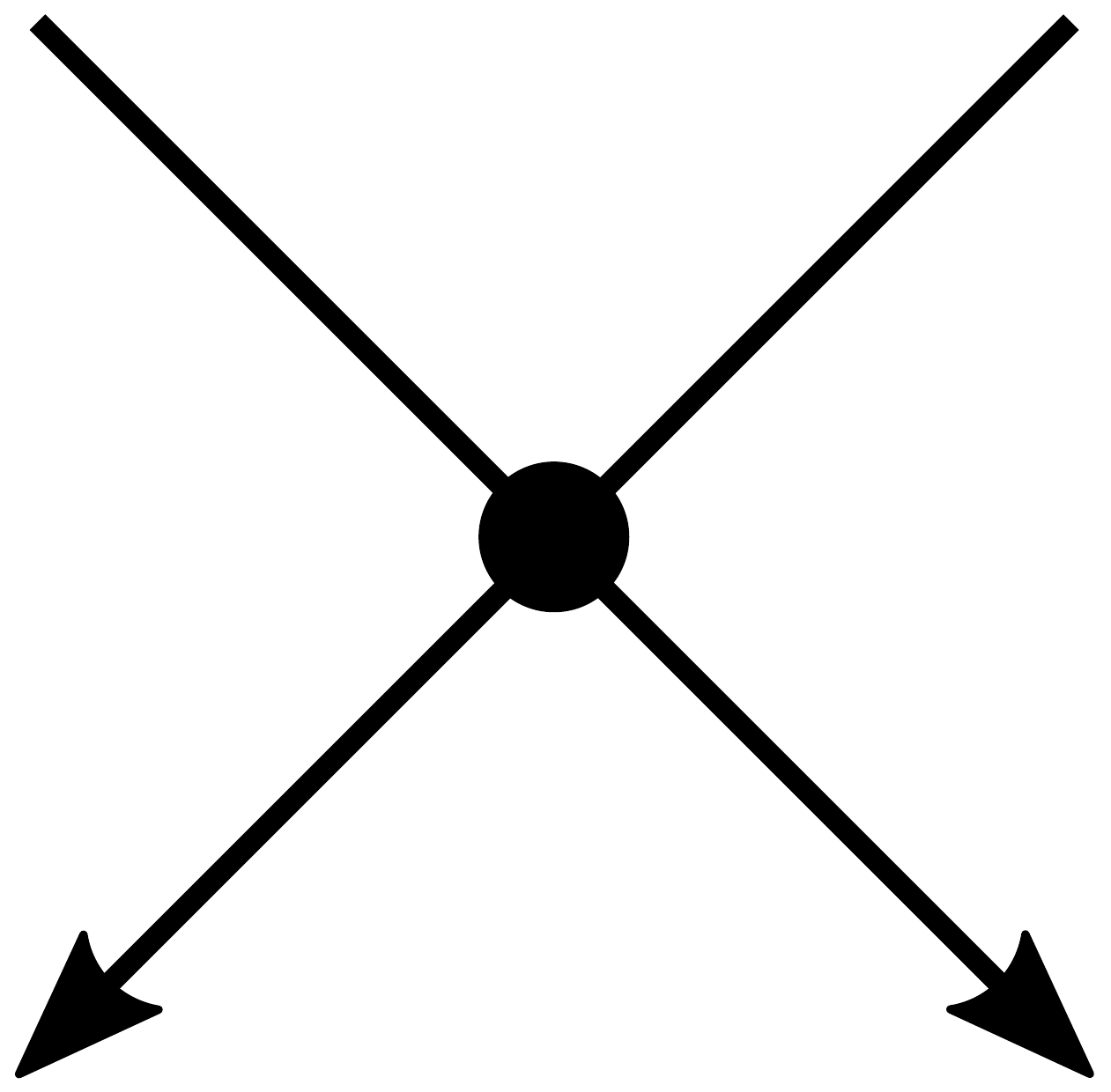}
\subcaption{Virtual crossing}
\end{center}
\end{minipage}
\end{tabular}
\caption{Classical and virtual crossings}
\label{crossing}
\end{figure}

\begin{defi}
A \textbf{virtual $n$-string link diagram} is an immersion of $n$ oriented intervals in $[0,1] \times [0,1]$ such that
\begin{itemize}
\item the initial point of each interval coincides with some $z_i \times \{ 0 \}$ and the endpoint coincides with $z_j \times \{ 1 \}$;
\item the singular set is a finite number of transverse double points;
\item each double point is labeled, either as a \textbf{positive crossing}, as a \textbf{negative crossing}, or as a \textbf{virtual crossing}; see Figure~\ref{crossing}.
\end{itemize}
The interval with initial point $z_i \times \{ 0 \}$ is called the \textbf{$i$-th string}.
The positive and negative crossings are also called the \textbf{classical crossings}.
A virtual string link diagram without the virtual crossing is said to be \textbf{classical}.
\end{defi}

\begin{figure}[h]
\begin{tabular}{ccc}
\begin{minipage}{0.33\hsize}
\begin{center}
\includegraphics[height=70pt]{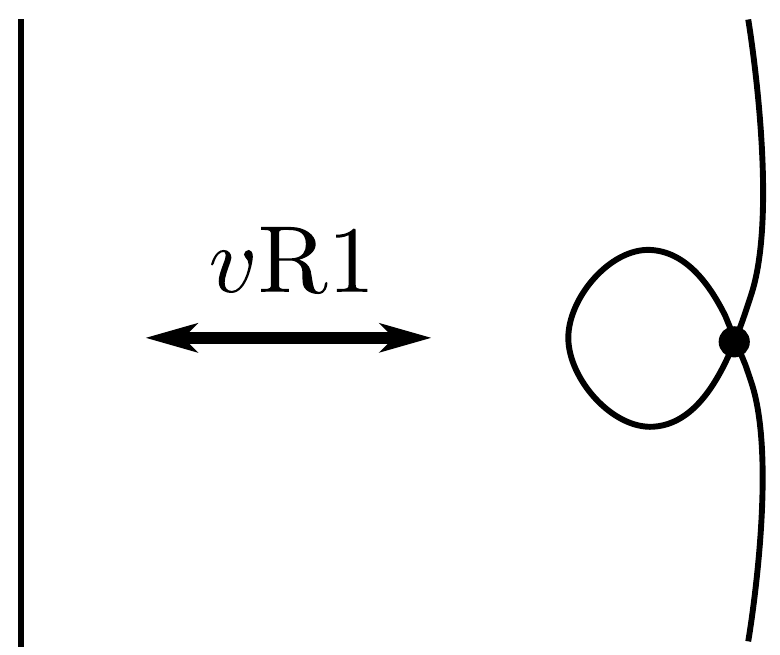}
\subcaption{}
\end{center}
\end{minipage}
\hspace{-20pt}
\begin{minipage}{0.33\hsize}
\begin{center}
\includegraphics[height=70pt]{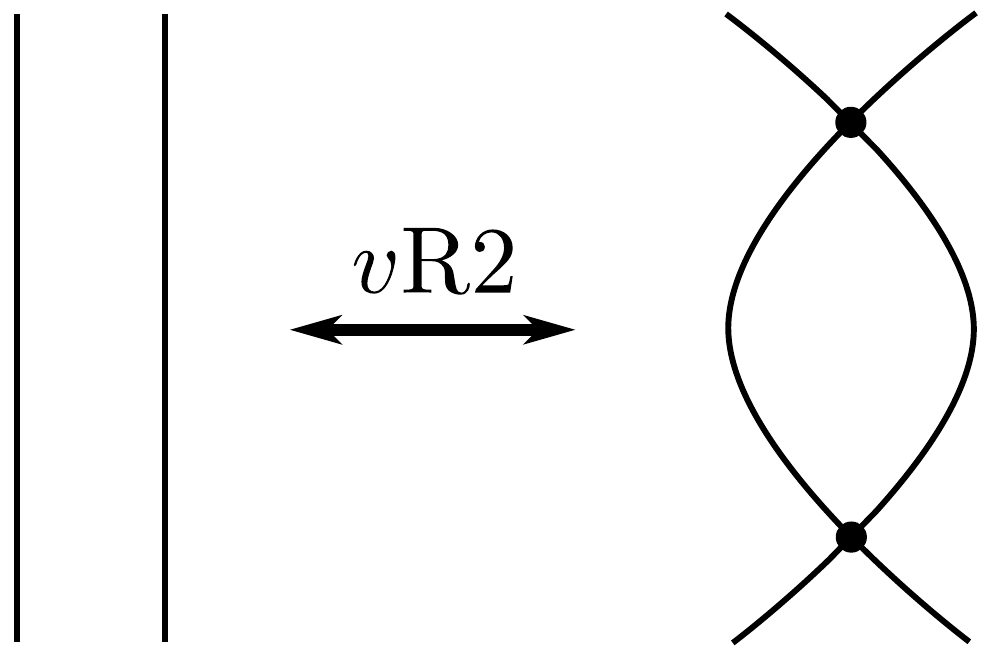}
\subcaption{}
\end{center}
\end{minipage}
\begin{minipage}{0.33\hsize}
\begin{center}
\includegraphics[height=70pt]{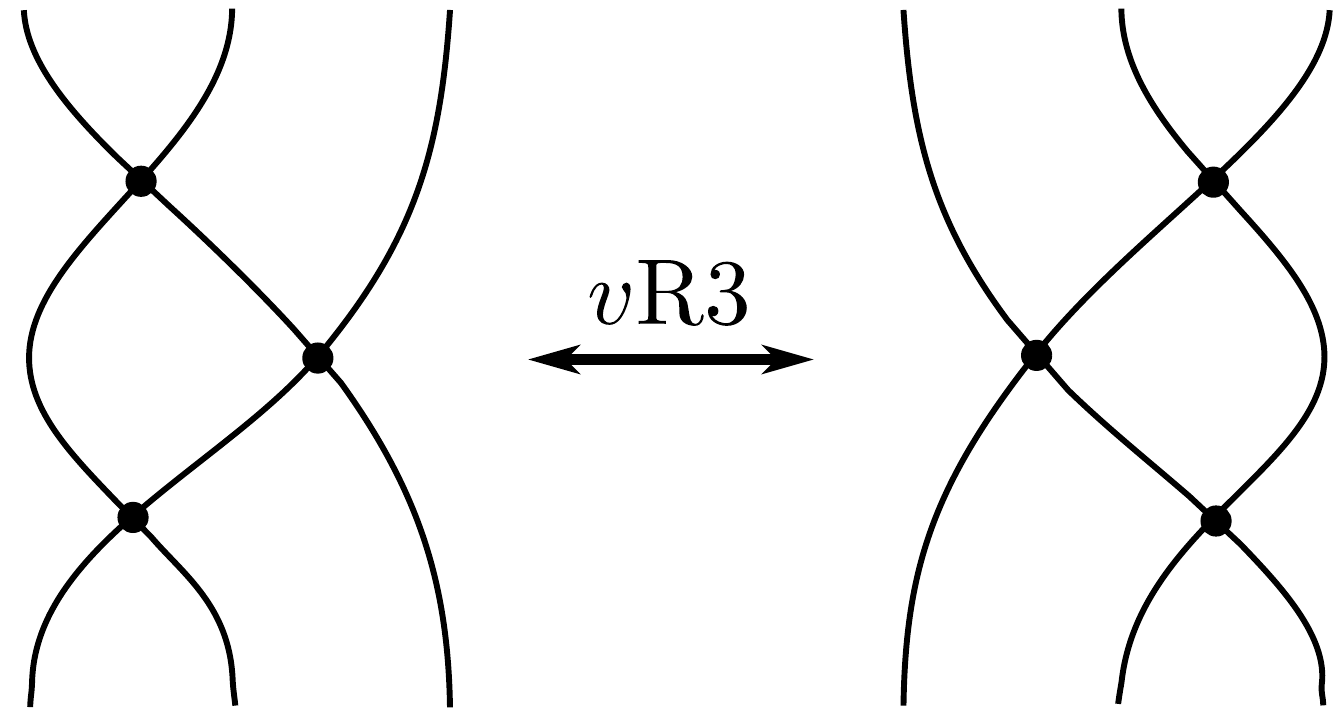}
\subcaption{}
\end{center}
\end{minipage}
\end{tabular}
\vspace{-15pt}
\caption{Virtual Reidemeister moves}
\label{virtual}
\end{figure}

\begin{figure}[h]
\hspace{20pt}
\begin{tabular}{cc}
\begin{minipage}{0.5\hsize}
\begin{center}
\includegraphics[height=70pt]{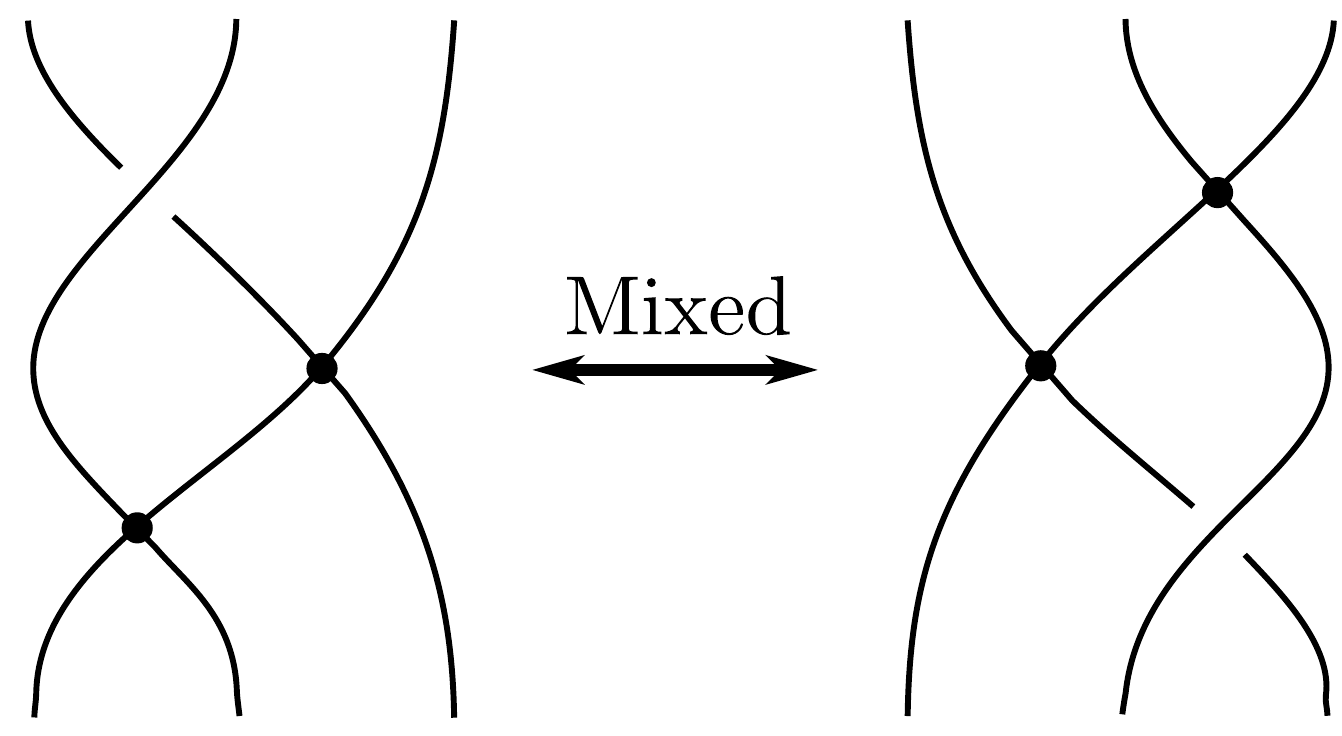}
\subcaption{}
\end{center}
\end{minipage}
\hspace{-50pt}
\begin{minipage}{0.5\hsize}
\begin{center}
\includegraphics[height=70pt]{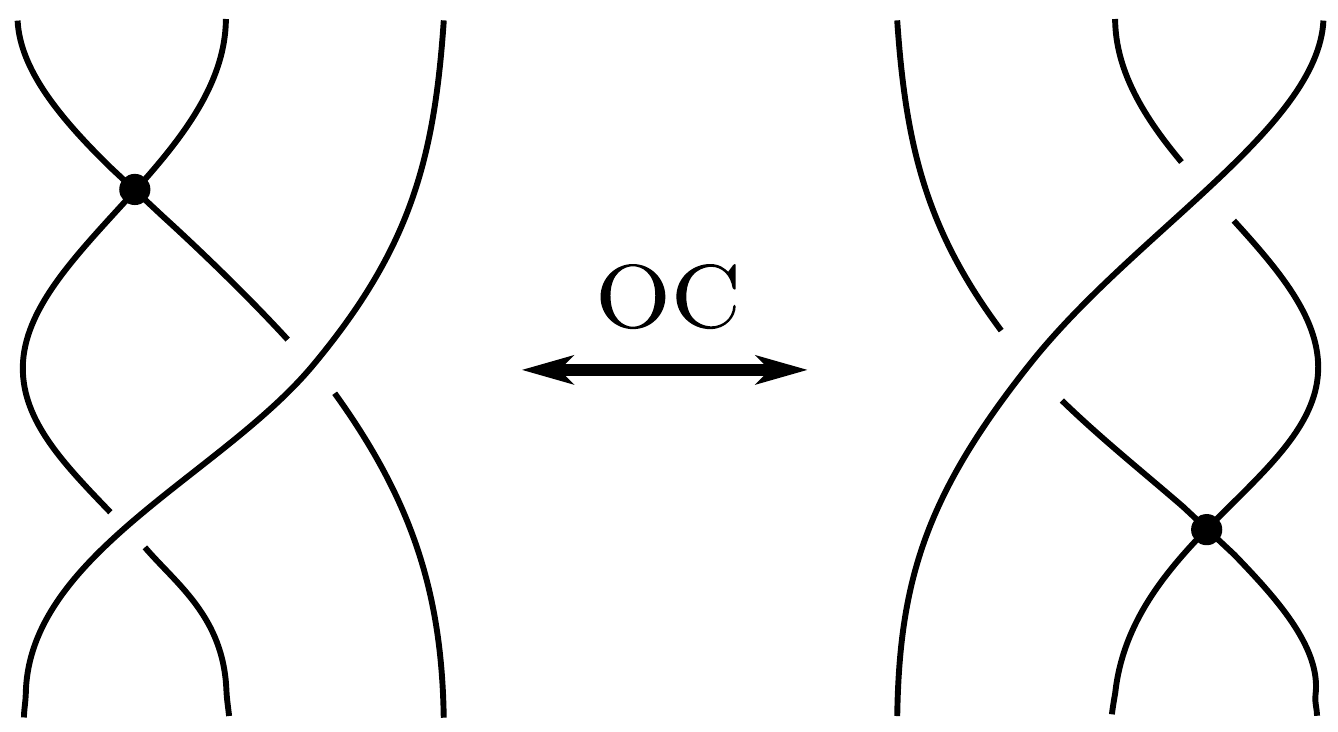}
\subcaption{}
\end{center}
\end{minipage}
\end{tabular}
\vspace{-15pt}
\caption{Mixed and OC moves}
\label{Mixed_OC}
\end{figure}

\begin{defi}
A \textbf{welded $n$-string link} is an equivalence class of the set of all virtual $n$-string link diagrams under the ambient isotopy on the plane and the welded Reidemeister moves, that is, the Reidemeister moves (see Figures~\ref{R1} and \ref{R2R3}), the virtual Reidemeister moves (see Figure~\ref{virtual}), the mixed move, and the over-commute (OC) move (see Figure~\ref{Mixed_OC}).
\end{defi}

The set of all welded $n$-string links has a structure of a monoid by concatenation and reparametrization, and we denote it by $w\mathcal{SL}_n$.
The $n$-string link monoid $\mathcal{SL}_n$ is naturally included in $w\mathcal{SL}_n$.

\begin{figure}[h]
\begin{center}
\includegraphics[height=100pt]{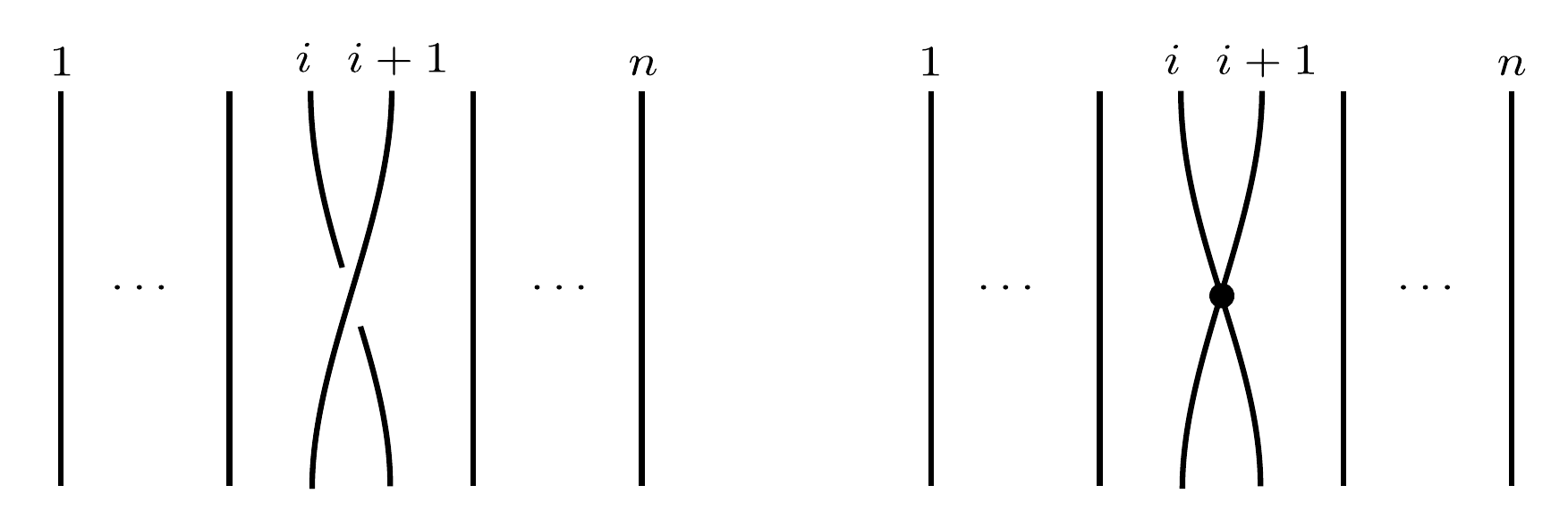}
\caption{Generators $\sigma_i$ and $\tau_i$}
\label{welded_diagram}
\end{center}
\end{figure}

If each string transversely meets every interval $[0,1] \times \{ t \}$ in a single point, the welded $n$-string link is also called the \textbf{welded $n$-braid}.
We denote $\WB_n$ by the set of all welded $n$-braid.
$\WB_n$ has the structure of a group, thus it is called the \textbf{welded braid group}.
The welded braid group has a presentation with generators $\{ \sigma_i,\tau_i \}_{i = 1,\ldots ,n-1}$ (see Figure~\ref{welded_diagram}) together with relations:
$$
\begin{cases}
\sigma_i \sigma_j = \sigma_j \sigma_i & (|i-j| \geq 2) \\ 
\sigma_i \sigma_{i+1} \sigma_i = \sigma_{i+1} \sigma_i \sigma_{i+1} & (i = 1, \ldots ,n-2) \\
\tau_i \tau_j = \tau_j \tau_i & (|i-j| \geq 2) \\
\tau_i \tau_{i+1} \tau_i = \tau_{i+1} \tau_{i} \tau_{i+1} & (i = 1, \ldots ,n-2) \\
\tau_i^2 = 1 &  (i = 1,\cdots ,n-1) \\
\sigma_i \tau_j = \tau_j \sigma_i & (|i-j| \geq 2) \\ 
\sigma_i \tau_{i+1} \tau_i = \tau_{i+1} \tau_i \sigma_{i+1} & (i = 1, \ldots ,n-2) \\
\tau_i \sigma_{i+1} \sigma_{i} = \sigma_{i+1} \sigma_{i} \tau_{i+1} & (i = 1, \ldots ,n-2).
\end{cases}
$$

\begin{remark}
The welded braid group $\WB_n$ is a generalization of the braid group.
Namely, the braid group $\B_n$ is interpreted as the fundamental group of the configuration space of $n$ points in the plane $\mathbb{R}^2$.
On the other hand, $\WB_n$ is interpreted as the fundamental group of the configuration space of $n$ Euclidean, unordered, disjoint, unlinked circles in the $3$-ball $B^3$ lying on planes parallel to a fixed one; see \cite{Damiani} for instance.
\end{remark}

The Tong-Yang-Ma representations for the braid groups may be extended to representations of the welded braid groups; see \cite{Bellingeri-Soulie} for example. Namely, it is defined for all $1\leq i\leq n-1$ by
$$
\sigma_i \longmapsto
I_{i-1} \oplus
\left(
\begin{array}{cc}
0 & 1 \\
t & 0
\end{array}
\right)
\oplus I_{n-i-1}
\ \ \textrm{and} \ \ 
\tau_i \longmapsto
I_{i-1} \oplus
\left(
\begin{array}{cc}
0 & \alpha^{-1} \\
\alpha & 0
\end{array}
\right)
\oplus I_{n-i-1}.
$$
We call it the \textbf{welded Tong-Yang-Ma representation}, and denote it by $wTYM_n \colon \WB_n \rightarrow \mathrm{GL}_{n} (\mathbb{Z} [t^{\pm1},\alpha^{\pm1}])$.

\subsection{Extension of welded Tong-Yang-Ma representation}
We recall from Section~\ref{ex} that given a diagram $D$ of a welded $n$-string link $L$, for each arc of $D$, we put vertex and label on it. Here, we regard the over-crossing and the virtual crossing arcs as the union of two and four arcs, respectively.
Moreover, assume that the labels of the top and bottom vertices are ordered just as the order of strings, and middle vertices are labeled by $m_1, \ldots, m_{2N-n}$.
Let $\mathcal{S}$ and $\mathcal{T}$ be the ordered $n$-tuples of top and bottom vertex labels, respectively.

\begin{figure}[h]
\begin{center}
\includegraphics[height=110pt]{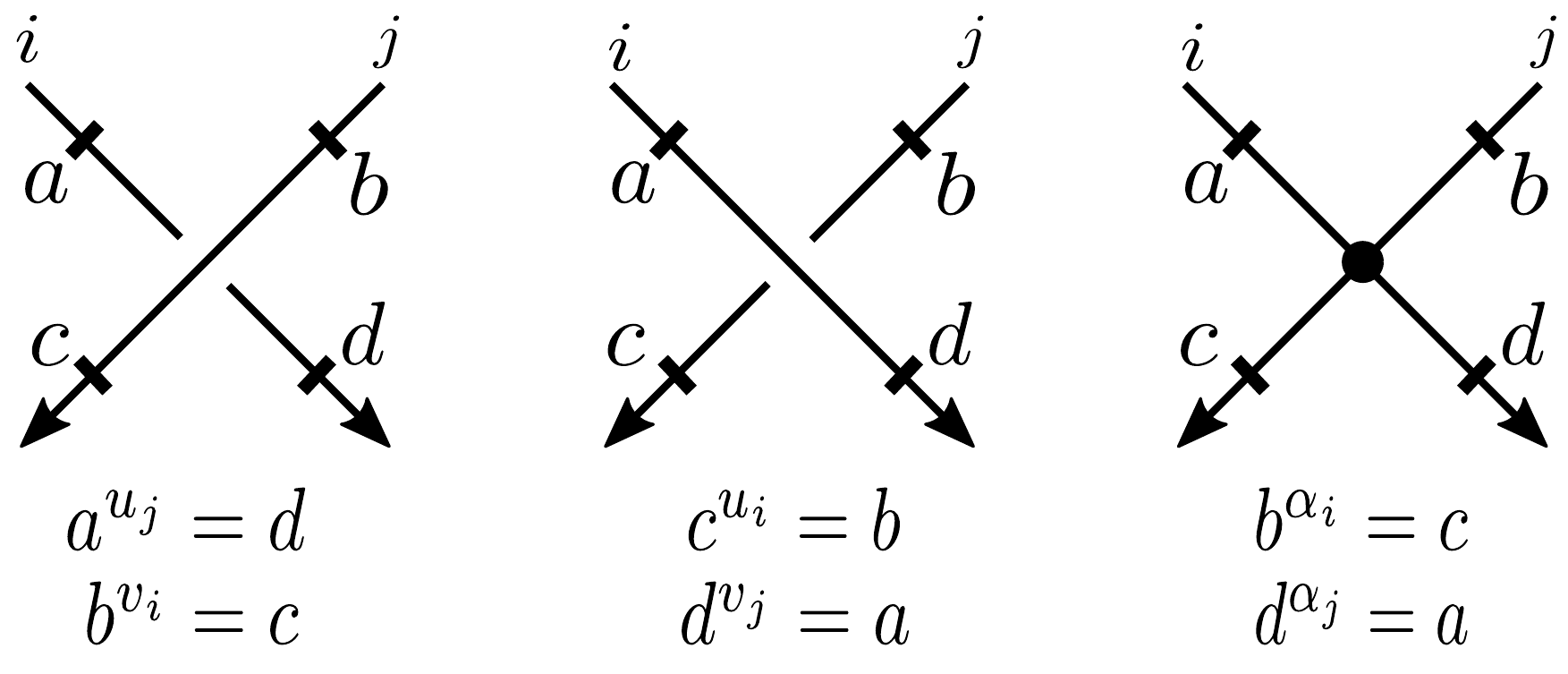}
\caption{Welded $\widetilde{\Lambda}_n$-relations}
\label{welded_Multi_TYM_crossing}
\end{center}
\end{figure}

Let $\widetilde{\Lambda}_n$ be the free abelian group of rank 3$n$ generated by $u_1, \ldots, u_n, v_1, \ldots, v_n, \alpha_1, \ldots, \alpha_n$.
Similarly to the classical case, we construct a $\widetilde{\Lambda}_n$-group associated to $D$ as follows.
\begin{defi}
Let $wK_n$ be the $\widetilde{\Lambda}_n$-group defined by the presentation:
\begin{description}
\item[Generators] the letters indexed by $\mathbb{Z}^{\oplus 3n}$ corresponding to each label of the vertex of $D$.
For example, for the label $a_1$, the set $\{ (a_1)_{\bm{n}} \}_{\bm{n} \in \mathbb{Z}^{\oplus 3n}}$ is the generating set.
We denote this set by the same letter $a_1$.
\item[Relations] the families of relations indexed by $\mathbb{Z}^{\oplus3n}$ of Figure~\ref{welded_Multi_TYM_crossing}.
\end{description}
\end{defi}

It is clear from its definition that this group is invariant under the welded Reidemeister moves except for R1.
In order for this last invariance to hold, we define the $\widetilde{\Lambda}_n$-group $w\widetilde{K}_n$ as follows. We may assume that $wK_n$ is of the form
$$wK_n = \left \langle a_1, \ldots, a_n,x_1, \ldots, x_n \relmiddle | x_1 = a_{i_1}^{w_1}, \ldots, x_n = a_{i_n}^{w_n} \right \rangle,$$
where $w_1, \ldots, w_n \in \widetilde{\Lambda}_n$ correspond to the welded $\widetilde{\Lambda}_n$-relations of $wK_n$.
Then, we define the $\widetilde{\Lambda}_n$-group $w\widetilde{K}_n$ with the following presentation
$$
w\widetilde{K}_n := \left \langle a_1, \ldots, a_n,x_1, \ldots, x_n \relmiddle | x_1 = a_{i_1}^{w_1 (u_{i_1}v_{i_1})^{-k_1}}, \ldots, x_n = a_{i_n}^{w_n (u_{i_n}v_{i_n})^{-k_n}} \right \rangle,
$$
where $k_j$ is the sum of signs of self-classical crossings in $i_j$-th string.

\begin{remark}
Similar works to define other groups of welded links have been carried out by Bardakov and Bellingeri \cite[Section 5]{BardakovBellingeri} and Bardakov \cite[Section 5]{Bardakov1}. However, as far as the we know, the $\widetilde{\Lambda}_n$-group $w\widetilde{K}_n$ does not appear in the literature.
\end{remark}

The $\widetilde{\Lambda}_n$-group $w\widetilde{K}_n$ is an invariant of a welded $n$-string link $L$, and its abelianization $w\widetilde{K}_n^{\textrm{ab}}$ is a free $\widetilde{\Lambda}_n$-module of rank $n$ with ordered bases $\mathcal{S}$ and $\mathcal{T}$.

\begin{defi}\label{def:welded multi-variable Tong-Yang-Ma matrix}
The \textbf{welded multi-variable Tong-Yang-Ma matrix} of a welded $n$-string link $L$ is defined as the change of basis matrix $w\mathcal{TYM}_n (L)$ such that $\mathcal{S}(w \mathcal{TYM}_n (L)) = \mathcal{T}$.
\end{defi}

The following result follows by repeating verbatim the proof of Theorem~\ref{mTYM},

\begin{thm}
\textup{(1)} The matrix $w\mathcal{TYM}_n (L)$ is invariant of $L$ for any welded $n$-string link $L$.

\textup{(2)} For all welded $n$-string links $L$ and $L'$, we have
$$
w\mathcal{TYM}_n (LL') = w\mathcal{TYM}_n (L) (L \cdot w\mathcal{TYM}_n (L')).
$$
\end{thm}

\begin{cor}
Let $w\mathcal{PSL}_n$ be the welded pure $n$-string link monoid.
Then, the map
$$
w\mathcal{TYM}_n \colon {\rm w} \mathcal{PSL}_n \longrightarrow \mathrm{GL}_n (\mathbb{Z}[\widetilde{\Lambda}_n]) {;}\ L \longmapsto w\mathcal{TYM}_n (L)
$$
is a representation of $w\mathcal{PSL}_n$.
\end{cor}

\begin{figure}[h]
\begin{center}
\includegraphics[height=100pt]{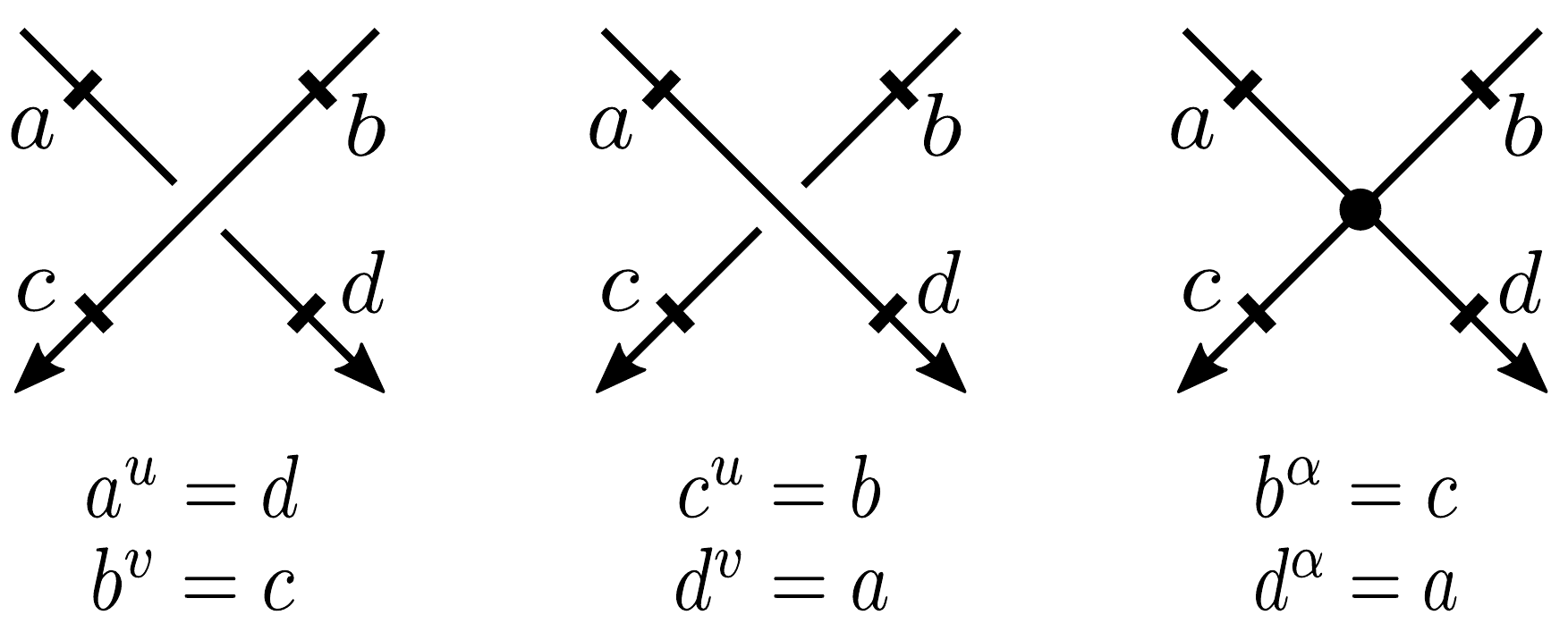}
\caption{Welded $\widetilde{\Lambda}$-relations}
\label{welded_TYM_crossing}
\end{center}
\end{figure}

If we replace $\widetilde{\Lambda}_n$ by the free abelian group $\widetilde{\Lambda}$ of rank 3 generated by $u,v$ and $\alpha$ and consider the welded $\widetilde{\Lambda}$-relations (see Figure~\ref{welded_TYM_crossing}), we also define the matrix $w\mathcal{TYM}(L)$ that we call the \textbf{welded three variable Tong-Yang-Ma matrix} in the same way as Section~\ref{two}.
Also, the specialization $u_1 = \cdots = u_n = u$, $v_1 = \cdots = v_n = v$ and $\alpha_1 = \cdots = \alpha_n = \alpha$ of $w\mathcal{TYM}_n (L)$ recovers $w\mathcal{TYM} (L)$.
It is not difficult to check that the map
$$
w \mathcal{TYM} \colon w \mathcal{SL}_n \longrightarrow \mathrm{GL}_n (\mathbb{Z}[\widetilde{\Lambda}]) {;}\ L \longmapsto w \mathcal{TYM} (L)
$$
is a representation of $w\mathcal{SL}_n$.

By straightforward computations, we deduce that:
\begin{align*}
&w \mathcal{TYM} (\sigma_i) = I_{i-1} \oplus
\left(
\begin{array}{cc}
0 & u \\
v & 0
\end{array}
\right)
\oplus I_{n-i-1}
\ \  \textrm{and}\\
&w \mathcal{TYM} (\tau_i) = I_{i-1} \oplus
\left(
\begin{array}{cc}
0 & \alpha^{-1} \\
\alpha & 0
\end{array}
\right)
\oplus I_{n-i-1}.
\end{align*}
Then, when we substitute $u = 1$ and $v = t$, these matrices coincide with the welded Tong-Yang-Ma representation.

\begin{figure}[h]
\begin{center}
\includegraphics[height=70pt]{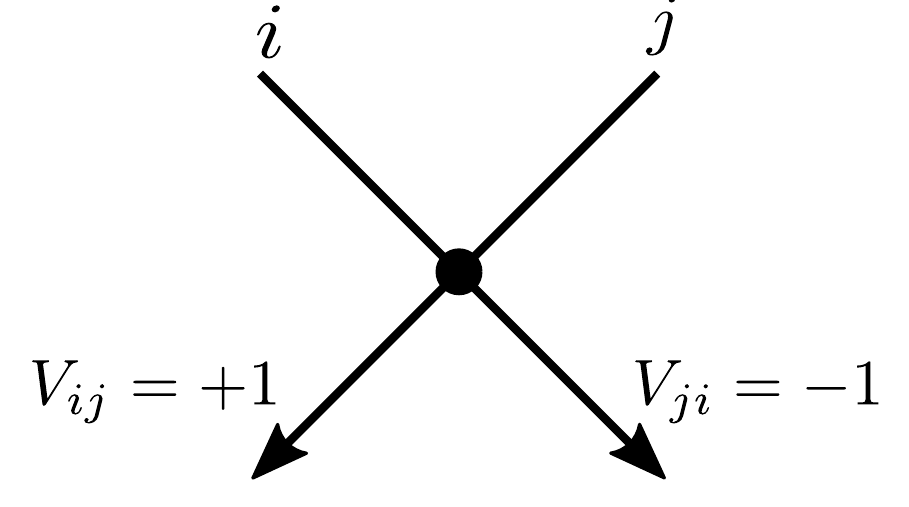}
\caption{The $i$-th string passes the $j$-th one ``from right to left'' and the $j$-th string passes the $i$-th one ``from left to right''.}
\label{virtualizing}
\end{center}
\end{figure}

We now determine the kernel of $w\mathcal{TYM} \colon w \mathcal{SL}_n \longrightarrow \mathrm{GL}_n (\mathbb{Z}[\widetilde{\Lambda}])$.
In order to do this, we introduce the linking number and virtual linking number of the welded string link.
In \cite{Goussarov-Polyak-Viro}, for a welded $n$-string link $L$, the \textbf{linking numbers} $v \ell_{ij} (L)$ and $v \ell_{ji} (L)$ are defined as follows: $v \ell_{ij} (L)$ is the sum of signs of the classical crossings where the $i$-th string of $L$ passes over the $j$-th one, while $v \ell_{ji} (L)$ is defined by exchanging the strings in the definition of $v \ell_{ij} (L)$.
Furthermore, we define the \textbf{virtual linking numbers} $V_{ij} (L)$ and $V_{ji} (L)$ in the following way: $V_{ij} (L)$ is given by subtracting the number of the virtual crossings where the $i$-th string passes the $j$-th string ``from left to right'' from the one ``from right to left'' (Figure~\ref{virtualizing}). Then, $V_{ji} (L)$ is defined by exchanging the strings in the definition of $V_{ij} (L)$, and we see that $V_{ij} (L) = -V_{ji} (L)$.

Using the same argument as the one in Section~\ref{ker}, we assume that the $\widetilde{\Lambda}$-group $w\widetilde{K}$ associated to $L$ has the following presentation:
$$
w \widetilde{K} = \left \langle a_1, \ldots, a_n,x_1, \ldots, x_n \relmiddle | x_1 = a_1^{u^{p_1} v^{q_1} \alpha^{r_1}}, \ldots, x_n = a_n^{u^{p_n} v^{q_n} \alpha^{r_n}} \right \rangle,
$$
where $p_i, q_i, r_i \in \mathbb{Z}$ for all $1\leq i\leq n$.
By the welded $\widetilde{\Lambda}$-relations, we can show that$$
p_j = \sum_{i = 1, i \neq j}^{n} v \ell_{ij} (L),
\ \ 
q_j = \sum_{i = 1, i \neq j}^{n} v \ell_{ji} (L)
\ \ \textrm{and}\ \ 
r_j = \sum_{i = 1, i \neq j}^{n} V_{ij} (L)
$$
for each $1\leq j \leq n$. Therefore, we deduce the following result.

\begin{thm}\label{thm:Ker_3}
For $n \geq 2$,
\begin{align*}
& \ker w\mathcal{TYM} 
\\
& = \left \{ L \in w\mathcal{SL}_n \relmiddle |  \mathrm{for \ each} \ j \ (j = 1, \ldots n), \ \sum_{i = 1, i \neq j}^{n} v \ell_{ij} (L) = \sum_{i = 1, i \neq j}^{n} v \ell_{ji} (L) = \sum_{i = 1, i \neq j}^{n} V_{ij} (L) = 0 \right \}.
\end{align*}
\end{thm}

Moreover, when we consider the multi-variable case, the kernel may be described as follows.

\begin{thm}\label{thm:Ker_4}
For $n \geq 2$,
$$
\ker w\mathcal{TYM}_n = \left \{ L \in w\mathcal{PSL}_n \relmiddle |  \mathrm{for \ any} \ 1 \leq i \neq j \leq n,\ v \ell_{ij} (L) = v \ell_{ji} (L) = V_{ij} (L) = 0 \right \}.
$$
\end{thm}

\subsection{Interpretation as biquandle}
In this subsection, we interpret the construction of the three-variable welded Tong-Yang-Ma matrix using the notion of a welded birack.
\begin{defi}
A \textbf{welded biquandle} is a non-empty set $X$ with four binary operations $\overline{\triangleleft}, \underline{\triangleleft}, \overline{\bullet}, \underline{\bullet} \colon X \times X \longrightarrow X$ satisfying the following axioms:
\begin{itemize}
\item[(0)] For any $y \in X$, the map $\overline{\triangleleft} y \colon X \rightarrow X$ sending $x$ to $x \overline{\triangleleft} y$ is bijective. \\
For any $y \in X$, the map $\underline{\triangleleft} y \colon X \rightarrow X$ sending $x$ to $x \underline{\triangleleft} y$ is bijective.
\item[(1)] For any $x \in X$,
$$
x \underline{\triangleleft}^{-1} x = x \overline{\triangleleft} (x \underline{\triangleleft}^{-1} x),
$$
where for any $y \in X$, $\underline{\triangleleft}^{-1} y$ (resp. $\overline{\triangleleft}^{-1} y$) $\colon X \rightarrow X$ is the inverse map of $\underline{\triangleleft} y$ (resp. $\overline{\triangleleft} y$).
\item[(2)] The map $S \colon X \times X \rightarrow X \times X$ defined by $S(x,y) := (y \overline{\triangleleft} x, x \underline{\triangleleft} y)$ is bijective.
\item[(3)] The map $S$ satisfies the set-theoretic Yang-Baxter equation:
$$
(S \times id_X) (id_X \times S) (S \times id_X) = (id_X \times S) (S \times id_X) (id_X \times S).
$$
\item[(v0)] For any $y \in X$, the map $\overline{\bullet} y \colon X \rightarrow X$ sending $x$ to $x \overline{\bullet} y$ is bijective. \\
For any $y \in X$, the map $\underline{\bullet} y \colon X \rightarrow X$ sending $x$ to $x \underline{\bullet} y$ is bijective.
\item[(v1)] For any $x \in X$,
$$
x \underline{\bullet}^{-1} x = x \overline{\bullet} (x \underline{\bullet}^{-1} x),
$$
where for any $y \in X$, $\underline{\bullet}^{-1} y$ (resp. $\overline{\bullet}^{-1} y$) $\colon X \rightarrow X$ is the inverse map of $\underline{\bullet} y$ (resp. $\overline{\bullet} y$).
\item[(v2)] The map $V \colon X \times X \rightarrow X \times X$ defined by $V(x,y) := (y \overline{\bullet} x, x \underline{\bullet} y)$ is bijective and satisfies $V^2 = id_{X\times X}$.
\item[(v3)] The map $V$ satisfies the set-theoretic Yang-Baxter equation:
$$
(V \times id_X) (id_X \times V) (V \times id_X) = (id_X \times V) (V \times id_X) (id_X \times V).
$$
\item[(M)] The maps $S$ and $V$ satisfy the following equality:
$$
(S \times id_X) (id_X \times V) (V \times id_X) = (id_X \times V) (V \times id_X) (id_X \times S).
$$
\item[(OC)] The maps $S$ and $V$ satisfy the following equality:
$$
(V \times id_X) (id_X \times S) (S \times id_X) = (id_X \times S) (S \times id_X) (id_X \times V).
$$
\end{itemize}
A \textbf{welded birack} is a non-empty set $X$ with four binary operations $\overline{\triangleleft}, \underline{\triangleleft}, \overline{\bullet}, \underline{\bullet} \colon X \times X \rightarrow X$ satisfying the above axioms except (1) and (v1).
\end{defi}

\begin{remark}
A \textbf{biquandle} (resp. \textbf{birack}) is a non-empty set $X$ with two binary operations $\overline{\triangleleft}, \underline{\triangleleft} \colon X \times X \rightarrow X$ satisfying axioms (0), (1), (2) and (3) (resp. (0), (2) and (3)).
Moreover, a \textbf{virtual biquandle} (resp. \textbf{virtual birack}) is a non-empty set $X$ with four binary operations $\overline{\triangleleft}, \underline{\triangleleft},
\overline{\bullet}, \underline{\bullet} \colon X \times X \rightarrow X$ satisfying the above axioms except (v5) (resp. (1), (v1) and (v5)).
\end{remark}

\begin{ex}
Let $M$ be a $\mathbb{Z}[u^{\pm1}, v^{\pm1}, \alpha^{\pm1}]$-module, then $x \overline{\triangleleft} y := ux, x \underline{\triangleleft} y := vx + (1-uv)y, x
\overline{\bullet} y := \alpha x$ and $x \underline{\bullet} y := \alpha^{-1} x$ define a virtual biquandle $(M, \overline{\triangleleft}, \underline{\triangleleft}, \overline{\bullet}, \underline{\bullet})$.
This is called the \textbf{(virtual) Alexander biquandle}.
The inverse map $S^{-1} \colon X \times X \rightarrow X \times X$ of $S$ is given by $S^{-1}(x,y) = (v^{-1}y + (1-u^{-1}v^{-1})x, u^{-1}x)$.

Moreover, if $x \overline{\triangleleft} y := vx, x \underline{\triangleleft} y := ux, x \overline{\bullet} y := \alpha x$ and $x \underline{\bullet} y := \alpha^{-1} x$, then
$(M, \overline{\triangleleft}, \underline{\triangleleft}, \overline{\bullet}, \underline{\bullet})$ is a welded birack.
The map $S^{-1} \colon X \times X \rightarrow X \times X$ is given by $S^{-1}(x,y) = (u^{-1}y, v^{-1}x)$.
We call it the \textbf{welded Tong-Yang-Ma birack}.
\end{ex}

\begin{figure}[h]
\begin{center}
\includegraphics[height=110pt]{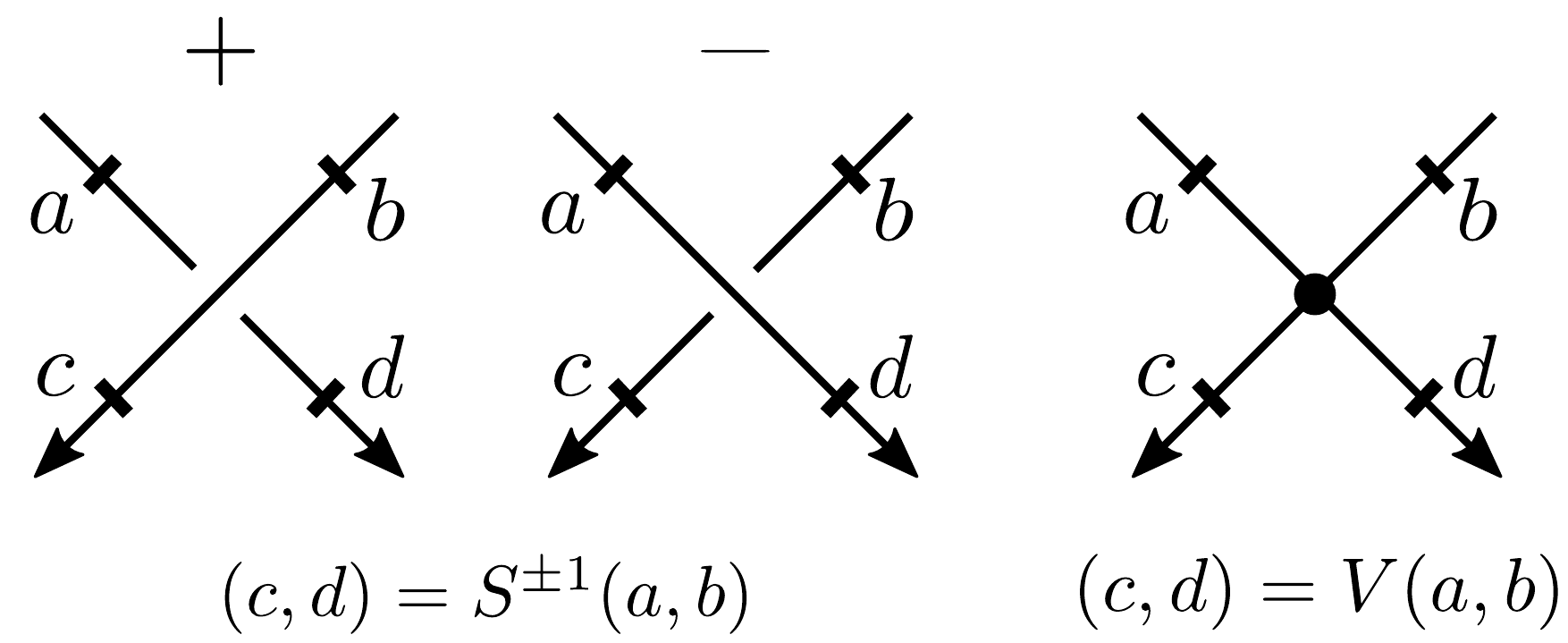}
\caption{welded biquandle relations}
\label{welded biquandle}
\end{center}
\end{figure}

Crans, Henrich and Nelson \cite{Crans-Henrich-Nelson} defined a virtual link invariant called the Alexander biquandle.
We construct a welded birack for a welded $n$-string link diagram by using their idea.
Let $D$ be a labeled welded $n$-string link diagram.
We define a welded birack $\mathcal{TYM}_{wbr}(D)$ associated to $D$ as follows:

\begin{defi} \label{def:welded Tong-Yang-Ma birack}
Let $\mathcal{TYM}_{wbr}(D)$ be a $\mathbb{Z}[u^{\pm1}, v^{\pm1}, \alpha^{\pm1}]$-module generated by labels of $D$ with the relations at positive ($+$) and negative ($-$) crossings and at virtual crossing described in Figure~\ref{welded biquandle}, where the maps $S, V \colon X \times X \rightarrow X \times X$ are given in the definition of the welded Tong-Yang-Ma birack.
\end{defi}

This module $\mathcal{TYM}_{wbr}(D)$ is just the abelianization of the $\widetilde{\Lambda}$-group associated to $D$.
Therefore, by correcting the module $\mathcal{TYM}_{wbr}(D)$ at self-classical crossings, the obtained module denoted by $\widetilde{\mathcal{TYM}}_{wbr}(D)$ is a free $\mathbb{Z}[u^{\pm1}, v^{\pm1}, \alpha^{\pm1}]$-module of rank $n$ and an invariant of a welded string link.
Moreover, the labels of the top and bottom vertices are both bases of $\widetilde{\mathcal{TYM}}_{wbr}(D)$, and the corresponding change of basis matrix is equal to the welded three-variable Tong-Yang-Ma matrix.

\section{Applying the Long-Moody construction to the Tong-Yang-Ma representation}\label{LM(TYM)}

Long and Moody \cite{Long} introduced a method to construct a new linear representation of the braid group $\mathbf{B}_{n}$ from a representation of $\mathbf{B}_{n+1}$.
This procedure complexifies the initial representation: for instance, it reconstructs the unreduced Burau representation from a one dimensional representation.
The underlying framework of this method was studied from a functorial point of view and extended in \cite{soulie1} and then generalised to other families of groups \cite{soulie2}.
In this section, we study the representation obtained from the Tong-Yang-Ma representations applying the Long-Moody construction. We fix a natural number $n\geq3$ all along Section~\ref{LM(TYM)}.

\subsection{The theoretical setting of the Long-Moody
construction}\label{sec:theoretical}

We detail here the required tool and present the abstract definition of a Long-Moody construction.

\subsubsection{Tool}
Recall that $\F_{n}=\left\langle x_{1},\ldots,x_{n}\right\rangle$ is the free group on $n$ generators.
We denote by $id_{1}\natural -$ the injection $\B_{n}\hookrightarrow\B_{n+1}$ defined by sending the generator $\sigma_i$ to $\sigma_{i+1}$ for each $1\leq i \leq n-1$.
The key ingredient to define the Long-Moody construction for braid groups is to find a group morphisms $\alpha_{n}\colon \B_{n}\longrightarrow\mathrm{Aut}(\F_{n})$ and $\chi_{n}\colon\F_{n}\longrightarrow\B_{n+1}$ such that:
\begin{itemize}
\item the morphism $\F_{n}*\B_{n}\longrightarrow\B_{n+1}$ given by the coproduct of $\chi_{n}$ and $id_{1}\natural -$ factors across the canonical surjection to the semidirect product $\F_{n}\rtimes_{\alpha_{n}}\B_{n}$;
\item the following diagram is commutative
\begin{equation}\label{eq:diagram}
\xymatrix{\F_{n}\ar@{^{(}->}[r]\ar@{->}[dr]_{\chi_{n}} &
\F_{n}\rtimes_{\alpha_{n}}\B_{n}\ar@{->}[d] &
\B_{n}\ar@{_{(}->}[l]\ar@{->}[dl]^{id_{1}\natural -}\\
 & \mathbf{B}_{n+1},
}
\end{equation}
where the vertical morphism $\F_{n}\rtimes_{\alpha_{n}}\B_{n}\longrightarrow\B_{n+1}$ is induced by the coproduct of $\chi_{n}$ and $id_{1}\natural -$.
\end{itemize}

In other words, we require that for all elements $\lambda\in \B_{n}$ and $x\in \F_{n}$ the morphism $\chi_{n}$ satisfies the following equality in $\B_{n+1}$:
\begin{equation}
(id_{1}\natural\lambda)\circ\chi_{n}(x)=\chi_{n}(\alpha_{n}(\lambda)(x))\circ(id_{1}\natural\lambda).\label{eq:cond1}
\end{equation}

\subsubsection{Definition}
A Long-Moody construction is defined as follows. We fix an abelian group $V$. We denote by $\mathcal{I}_{\F_{n}}$ the augmentation ideal of the group ring $\Z\left[\F_{n}\right]$.
We note that the action $\alpha_{n}$ canonically induces an action of $\B_{n}$ on $\mathcal{I}_{\F_{n}}$ (that we denote in the same way for convenience).

Let $\rho\colon\B_{n+1}\longrightarrow \mathrm{Aut}_{\Z}(V)$ be a linear representation.
The module $V$ has a canonical $\F_{n}$-module structure induced by precomposing the representation $\rho$ by the morphism $\chi_{n}$.
Then the Long-Moody construction defined by the morphisms $\alpha_{n}$ and $\chi_{n}$
$$
\mathbf{LM}_{\alpha_{n},\chi_{n}}(\rho):\B_{n}\longrightarrow
\mathrm{Aut}_{\Z}(\mathcal{I}_{\F_{n}}\underset{\F_{n}}{\otimes}V)
$$
is the map defined by:
$$
\mathbf{LM}_{\alpha_{n},\chi_{n}}(\rho)(\lambda)(i\underset{\F_{n}}{\otimes}v)=(\alpha_{n}(\lambda)(i)\underset{\F_{n}}{\otimes}\rho(id_{1}\natural\lambda)(v))
$$
for all $\lambda\in \B_{n}$, $i\in\mathcal{I}_{\F_{n}}$ and $v\in V$. For sake of completeness, we detail that:

\begin{lemma}[{\cite[Section 2.2.4]{soulie2}}]
The representation $\mathbf{LM}_{\alpha_{n},\chi_{n}}(\rho)$ is well-defined.
\end{lemma}

\begin{proof}
We consider elements $\lambda\in\B_{n}$, $x\in\F_{n}$, $v\in V$ and $i\in\mathcal{I}_{\F_{n}}$.
First, since $\rho$ is a morphism, we deduce from \eqref{eq:cond1} that
$$
\mathbf{LM}_{\alpha_{n},\chi_{n}}(\rho)(\lambda)(i\underset{\F_{n}}{\otimes}\rho(\chi_{n}(x))(v))=\mathbf{LM}_{\alpha_{n},\chi_{n}}(\rho)(\lambda)(i\cdot x\underset{\F_{n}}{\otimes}v),
$$
which gives the compatibility of the assignment $\mathbf{LM}_{\alpha_{n},\chi_{n}}(\rho)$ with respect to the tensor product over $\Z\left[\mathbf{F}_{n}\right]$. That this assignment $\mathbf{LM}_{\alpha_{n},\chi_{n}}(\rho)$ defines a morphism on $\B_{n}$ follows from the fact that $\alpha_{n}$ and $\rho$ are themselves morphisms.
\end{proof}

The original Long-Moody construction of \cite{Long} uses the Artin homomorphism $a_{n} \colon \B_n \longrightarrow \textrm{Aut} (\F_n)$ for the morphism $\alpha_{n}$:
$$
a_{n} (\sigma_{i}) (x_{j}) =
\begin{cases}
x_{i+1} & \textrm{if \ensuremath{j = i},}\\
x_{i+1}^{-1} x_{i} x_{i+1} & \textrm{if \ensuremath{j = i+1},}\\
x_{j} & \textrm{if \ensuremath{j\notin\{i, i+1\}},}
\end{cases}
$$
and we fix the assignment $\alpha_{n}=a_{n}$ from now on.
We could use another Wada representation for $a_{n}$; see \cite{Wada,Ito}.
We may always choose the trivial morphism $\F_{n}\rightarrow 0\rightarrow \B_{n+1}$ as $\chi_{n}$ for the relation \eqref{eq:cond1} to be satisfied.
However the construction is much more interesting using a non-trivial morphism for this parameter.
Indeed applying the Long-Moody construction with the trivial $\chi_{n}$ to a one-dimensional representation provides the permutation representation of $\B_{n}$ (sending the braid generators on the permutation matrix) and the iteration of this Long-Moody construction gives the tensor powers of that permutation representation; see \cite[Section 2.2.5]{soulie2} for further details.
The first non-trivial instance of such $\chi_{n}$ is the one used to define the original Long-Moody construction: let $\chi_{n,1}:\ensuremath{\F_{n}}\hookrightarrow\B_{n+1}$ be the injective morphism defined by
$$
x_{i}\longmapsto\left(\sigma_{i-1}\cdots\sigma_{2}\sigma_{1}\right)^{-1}\sigma_{i}^{2}\left(\sigma_{i-1}\cdots\sigma_{2}\sigma_{1}\right).
$$
The morphisms $a_{n}$ and $\chi_{n,1}$ satisfy the equality \eqref{eq:cond1}; see \cite[Proposition 2.8]{soulie1}.
They define the original Long-Moody construction of \cite{Long} that we denote by $\mathbf{LM}$ for simplicity.

\subsection{The application to the Tong-Yang-Ma representation}

We now study the application of the Long-Moody construction $\mathbf{LM}$ on the Tong-Yang-Ma representation $TYM_{n+1}$ recalled in Section~\ref{sec:tym}. More precisely, we consider the following slight modification of the Tong-Yang-Ma representation.
We generically denote by $\Z[r^{\pm1}]_{r}$ the braid group representation defined by the module $\Z[r^{\pm1}]$, where the action is induced by sending each Artin generator to the multiplication by $r$.
Also, the Tong-Yang-Ma representation can be viewed as a representation over the ring of polynomials with two variables $\Z[t^{\pm1},q^{\pm1}]$ by using the canonical inclusion $\Z[t^{\pm1}] \hookrightarrow \Z [t^{\pm1},q^{\pm1}]$, and thus consider $TYM_{n+1} \colon \B_{n+1} \longrightarrow \mathrm{GL}_{n+1} (\Z[t^{\pm1},q^{\pm1}])$ (we keep the same notation for sake of simplicity).
Then we actually apply the construction $\mathbf{LM}$ to the tensor product representation $\Z[q^{\pm1}]_{q}\otimes_{\Z[t^{\pm1},q^{\pm1}]}TYM_{n+1}$.

\subsubsection{General property of the Long-Moody construction}
Before describing the kernel of the application of the Long-Moody construction to the Tong-Yang-Ma representation, we need to introduce the following technical concepts.
First, we note that the representation $TYM_{n+1}$ may be seen as a $\B_{n}$-module using the canonical inclusion $id_{1}\natural -\colon \B_{n} \hookrightarrow \B_{n+1}$.
We denote by $\tau_{1}TYM_{n}$ that representation:
$$
\tau_{1}TYM_{n} \colon \B_{n} \hookrightarrow \B_{n+1} \longrightarrow \mathrm{GL}_{n+1} (\mathbb{Z} [t^{\pm1},q^{\pm1}]).
$$
In the same way, we denote by $\tau_{2}TYM_{n}$ the representation $TYM_{n+2}$ seen as a $\B_{n}$-module by precomposing by the canonical inclusion $(id_{1}\natural -) \circ (id_{1}\natural -)\colon \B_{n} \hookrightarrow \B_{n+1}\hookrightarrow \B_{n+2}$.
We also consider the inclusion $\Z[t^{\pm1},q^{\pm1}]^{\oplus n}\hookrightarrow\Z[t^{\pm1},q^{\pm1}]^{\oplus n+1}$ defined by sending the left-hand term on the $n$ last copies of $\Z[t^{\pm1},q^{\pm1}]$ of the right-hand term: it directly follows from the definition of the Tong-Yang-Ma representation that this induces a $\B_{n}$-module inclusion $i_{1}\colon TYM_{n} \hookrightarrow\tau_{1}TYM_{n}$ which kernel is isomorphic to the trivial representation $\Z[t^{\pm1},q^{\pm1}]_{1}$.

For sake of simplicity, we denote by $\mathbf{LM}_{q}(TYM_{n+1})$ the application of the Tong-Yang-Ma representation $\Z[q^{\pm1}]_{q^{-1}}\otimes_{\Z[t^{\pm1},q^{\pm1}]}\mathbf{LM}(\Z[q^{\pm1}]_{q}\otimes_{\Z[t^{\pm1},q^{\pm1}]}TYM_{n+1})$.
A key result of \cite{soulie1} is that Long-Moody construction applied to the morphism $i_{1}$ induces a morphism $\tilde{i}_{1}\colon\mathbf{LM}_{q}(TYM_{n+1})\longrightarrow \tau_{1}\mathbf{LM}_{q}(TYM_{n+1})$; see \cite[Sections 2.2]{soulie1}. 
Moreover, it follows from \cite[Theorem 4.23]{soulie1} that $\tilde{i}_{1}$ is injective and its cokernel is isomorphic to
$$
\tau_{2}TYM_{n}\oplus \mathbf{LM}_{q}(\Z[t^{\pm1}]_{1})
$$
as a $\B_{n}$-module. We refer the reader to \cite[Sections 2.2 and 4]{soulie1} for the detailed proofs of these properties.

We recall that $Bur_{n,t}$ denotes the unreduced Burau representation of the braid group $\B_{n}$ (here the extra index $t$ denotes the defining parameter of the ring of this representation); see Section~\ref{sec:tym}.
Then, we actually know that $\mathbf{LM}_{q}(\Z[t^{\pm1}]_{1})$ is isomorphic to $Bur_{n,q^{2}}$ as a $\B_{n}$-module; see \cite[Proposition 2.30]{soulie1}.
Hence, there is a $\B_{n}$-module surjection
\begin{equation}\label{eq:canonical_surjection}
\tau_{1}\mathbf{LM}_{q}(TYM_{n+1})\twoheadrightarrow Bur_{n,q^{2}}.
\end{equation}
Then we straightforwardly deduce from that surjection \eqref{eq:canonical_surjection} the following inclusion for the kernels of the representations:
\begin{thm}
There is a natural inclusion
$\mathrm{Ker}(\tau_{1}\mathbf{LM}_{q}(TYM_{n+1}))\subseteq
\mathrm{Ker}(Bur_{n,q^{2}})$.
\end{thm}

\subsubsection{Computations on the elements of the kernel of the
Burau representation}
We now show that the known elements of the kernel of the Burau representation (see \cite{Long-Paton} and \cite{Bigelow}) are not contained in the kernel of $\tau_{1}\mathbf{LM}_q(TYM_{n+1})$.

$\mathbf{For}$ $\boldsymbol{n=5:}$ we set $\psi_{1} :=\sigma_{3}^{-1} \sigma_{2} \sigma_{1}^{2} \sigma_{2}\sigma_{4}^{3} \sigma_{3} \sigma_{2}$, $\psi_{2} :=\sigma_{4}^{-1} \sigma_{3} \sigma_{2} \sigma_{1}^{-2} \sigma_{2}\sigma_{1}^{2} \sigma_{2}^{2} \sigma_{1} \sigma_{4}^{5}$ and
$$
\sigma := [\psi_{1}^{-1} \sigma_{4} \psi_{1}, \psi_{2}^{-1} \sigma_{4} \sigma_{3} \sigma_{2} \sigma_{1}^{2} \sigma_{2} \sigma_{3} \sigma_{4} \psi_{2}].
$$
Then $\sigma$ is a non-trivial element of $\ker (Bur_5)$, where $[a, b] := a^{-1} b^{-1} a b$ for $a, b \in \B_n$.

$\mathbf{For}$ $\boldsymbol{n=6:}$ we set $\delta_1 := \sigma_4\sigma_5^{-1} \sigma_2^{-1} \sigma_1$, $\delta_2 :=\sigma_4^{-1} \sigma_5^2 \sigma_2 \sigma_1^{-2}$ and $\theta := \sigma_{5}^{-1} \sigma_{4}^{-1} \sigma_{5} \sigma_{3}^{-1} \sigma_{4} \sigma_{2}^{-1} \sigma_{3}^{-3} \sigma_{1}^3 \sigma_{5} \sigma_{4} \sigma_{3}^{-1} \sigma_{2}^{-1} \sigma_{1}$.
Now we define
$$
\tau := [\delta_{1}^{-1} \sigma_{3} \delta_{1}, \delta_{2}^{-1} \sigma_{3} \delta_2]\ \textrm{and}\ \xi :=[\theta^{-1} \sigma_5 \theta, (\sigma_2 \sigma_3 \sigma_4 \sigma_5)^5].
$$
Then $\tau$ and $\xi$ are both non-trivial and in $\ker (Bur_6)$.

$\mathbf{For}$ $\boldsymbol{n=7:}$ we set $\upsilon := [\sigma_1, \gamma_6(\theta)^{-1} \sigma_6 \gamma_6(\theta)]$.
Then $\upsilon$ is also non-trivial and in $\ker (Bur_7)$.

Using Mathematica\footnote{See \url{https://www.dropbox.com/s/hyunrjx0zd0cnv0/Mathematica.zip?dl=0}}, we have checked that $\sigma$, $\tau$, $\xi$ and $\upsilon$ are elements of $\ker(\mathbf{LM}(TYM_{n+1}))$, but not of $\ker(\tau_{1}\mathbf{LM}_q(TYM_{n+1}))$ for $n=5,6$ and $7$.
This proves that:
\begin{cor}\label{coro:kernelLM}
There is a strict inclusion $\ker(\tau_{1}\mathbf{LM}_{q}(TYM_{n+1}))\subsetneq \ker(Bur_{n,q^{2}})$.
\end{cor}
Therefore, these observations raise the following natural question about these representations:
\begin{q} \label{kernel}
Is it true that $\ker (\mathbf{LM}(TYM_{n+1})) = \ker (Bur_n)$ and that $\ker(\tau_{1}\mathbf{LM}_{q}(TYM_{n+1}))$ is trivial ?
\end{q}

\subsubsection{A new family of representations of the braid groups}
We finally study the application of the Long-Moody construction $\mathbf{LM}_{q}(TYM_{n+1})$.
In particular, it induces a family of representations of the braid groups which appears to be new to us; see Proposition~\ref{prop:not_equiv_LKB}.
First of all, we reduce the representation $\mathbf{LM}_{q}(TYM_{n+1})$ as follows.
\begin{align*}
\mathbf{LM}_q(TYM_4)(\sigma_1) &=
\left(
\begin{array}{cccccccccccc}
 0 & 0 & 0 & 0 & q^2 t & 0 & 0 & 0 & 0 & 0 & 0 & 0 \\
 0 & 0 & 0 & 0 & 0 & 0 & q^2 & 0 & 0 & 0 & 0 & 0 \\
 0 & 0 & 0 & 0 & 0 & q^2 t^2 & 0 & 0 & 0 & 0 & 0 & 0 \\
 0 & 0 & 0 & 0 & 0 & 0 & 0 & q^2 & 0 & 0 & 0 & 0 \\
 1 & 0 & 0 & 0 & 1-q^2 t & 0 & 0 & 0 & 0 & 0 & 0 & 0 \\
 0 & 0 & 1 & 0 & 0 & 0 & 1-q^2 t & 0 & 0 & 0 & 0 & 0 \\
 0 & t & 0 & 0 & 0 & \left(1-q^2\right) t & 0 & 0 & 0 & 0 & 0 & 0 \\
 0 & 0 & 0 & 1 & 0 & 0 & 0 & 1-q^2 & 0 & 0 & 0 & 0 \\
 0 & 0 & 0 & 0 & 0 & 0 & 0 & 0 & 1 & 0 & 0 & 0 \\
 0 & 0 & 0 & 0 & 0 & 0 & 0 & 0 & 0 & 0 & 1 & 0 \\
 0 & 0 & 0 & 0 & 0 & 0 & 0 & 0 & 0 & t & 0 & 0 \\
 0 & 0 & 0 & 0 & 0 & 0 & 0 & 0 & 0 & 0 & 0 & 1 \\
\end{array}
\right),\\
\mathbf{LM}_q(TYM_4)(\sigma_2) &= \mathbf{LM}_q(TYM_4)(\sigma_1)^s,
\end{align*}
where $s$ is the permutation matrix for $(1,5,9) (2,7,12) (3,8,10) (4,6,11) \in \mathfrak{S}_{12}$ and $A^s$ denotes the conjugation of a square matrix $A$ by $s$.
Gathering together the first, $5$-th, and $9$-th rows and columns of these matrices, we transform them by a change of basis into the following block matrices: 
\begin{align*}
\mathbf{LM}_q(TYM_4)(\sigma_1) &=
\left(
\begin{array}{cccccccccccc}
\color{red} 0 & \color{red} q^2 t & \color{red} 0 & 0 & 0 & 0 & 0 & 0 & 0 & 0 & 0 & 0 \\
\color{red} 1 & \color{red} 1-q^2 t & \color{red} 0 & 0 & 0 & 0 & 0 & 0 & 0 & 0 & 0 & 0 \\
\color{red} 0 & \color{red} 0 & \color{red} 1 & 0 & 0 & 0 & 0 & 0 & 0 & 0 & 0 & 0 \\
 0 & 0 & 0 & \color{blue} 0 & \color{blue} 0 & \color{blue} 0 & \color{blue} 0 & \color{blue} q^2 & \color{blue} 0 & \color{blue} 0 & \color{blue} 0 & \color{blue} 0 \\
 0 & 0 & 0 & \color{blue} 0 & \color{blue} 0 & \color{blue} 0 & \color{blue} q^2 t^2 & \color{blue} 0 & \color{blue} 0 & \color{blue} 0 & \color{blue} 0 & \color{blue} 0 \\
 0 & 0 & 0 & \color{blue} 0 & \color{blue} 0 & \color{blue} 0 & \color{blue} 0 & \color{blue} 0 & \color{blue} q^2 & \color{blue} 0 & \color{blue} 0 & \color{blue} 0 \\
 0 & 0 & 0 & \color{blue} 0 & \color{blue} 1 & \color{blue} 0 & \color{blue} 0 & \color{blue} 1-q^2 t & \color{blue} 0 & \color{blue} 0 & \color{blue} 0 & \color{blue} 0 \\
 0 & 0 & 0 & \color{blue} t & \color{blue} 0 & \color{blue} 0 & \color{blue} \left(1-q^2\right) t & \color{blue} 0 & \color{blue} 0 & \color{blue} 0 & \color{blue} 0 & \color{blue} 0 \\
 0 & 0 & 0 & \color{blue} 0 & \color{blue} 0 & \color{blue} 1 & \color{blue} 0 & \color{blue} 0 & \color{blue} 1-q^2 & \color{blue} 0 & \color{blue} 0 & \color{blue} 0 \\
 0 & 0 & 0 & \color{blue} 0 & \color{blue} 0 & \color{blue} 0 & \color{blue} 0 & \color{blue} 0 & \color{blue} 0 & \color{blue} 0 & \color{blue} 1 & \color{blue} 0 \\
 0 & 0 & 0 & \color{blue} 0 & \color{blue} 0 & \color{blue} 0 & \color{blue} 0 & \color{blue} 0 & \color{blue} 0 & \color{blue} t & \color{blue} 0 & \color{blue} 0 \\
 0 & 0 & 0 & \color{blue} 0 & \color{blue} 0 & \color{blue} 0 & \color{blue} 0 & \color{blue} 0 &\color{blue} 0 & \color{blue} 0 & \color{blue} 0 & \color{blue} 1 \\
\end{array}
\right),\\
\mathbf{LM}_q(TYM_4)(\sigma_2) &= \mathbf{LM}_q(TYM_4)(\sigma_1)^u,
\end{align*}
where $u$ is the permutation matrix for $(1,2,3) (2,8,12) (5,9,10) (6,7,11)$.

Now, we consider the representation $\eta_n \colon \F_n \rtimes_{a_n} \B_n \longrightarrow \mathrm{GL}_n(\mathbb{Z}[t^{\pm1}, q^{\pm1}]))$ defined by
\begin{align*}
\sigma_i \longmapsto TYM_n(\sigma_i) \hspace{10pt} \textrm{and} \hspace{10pt} x_i \longmapsto q \Diag (1,\ldots, 1, \overset{i}{t}, 1, \ldots, 1).
\end{align*}
For $n=3$, the application of the Long-Moody construction $\mathbf{LM}_q(\eta_3) \colon \B_3 \longrightarrow \mathrm{GL}_9(\mathbb{Z}[t^{\pm1}, q^{\pm1}]))$ is defined by
\begin{align*}
\mathbf{LM}_q(\eta_3) (\sigma_1) &= \left(
\begin{array}{ccccccccc}
 0 & 0 & 0 & 0 & q^2 & 0 & 0 & 0 & 0 \\
 0 & 0 & 0 & q^2 t^2 & 0 & 0 & 0 & 0 & 0 \\
 0 & 0 & 0 & 0 & 0 & q^2 & 0 & 0 & 0 \\
 0 & 1 & 0 & 0 & 1-q^2 t & 0 & 0 & 0 & 0 \\
 t & 0 & 0 & \left(1-q^2\right) t & 0 & 0 & 0 & 0 & 0 \\
 0 & 0 & 1 & 0 & 0 & 1-q^2 & 0 & 0 & 0 \\
 0 & 0 & 0 & 0 & 0 & 0 & 0 & 1 & 0 \\
 0 & 0 & 0 & 0 & 0 & 0 & t & 0 & 0 \\
 0 & 0 & 0 & 0 & 0 & 0 & 0 & 0 & 1 \\
\end{array}
\right),\\
\mathbf{LM}_q(\eta_3) (\sigma_2) &= \mathbf{LM}_q(\eta_3) (\sigma_1)^v,
\end{align*}
where $v$ is the permutation matrix for $(1,5,9)(2,6,7)(3,4,8) \in \mathfrak{S}_9$.
These matrices are exactly the second diagonal block of $\mathbf{LM}_q(TYM_4)$.
Therefore, we obtain the decomposition
$$
\mathbf{LM}_q(TYM_4) \cong Bur_{3, q^2t} \oplus \mathbf{LM}_q(\eta_3),
$$
In general, the matrix presentation of $\mathbf{LM}_q(TYM_{n+1}) (\sigma_i)$ is given as follows:
\begin{align*}
\mathbf{LM}_q(TYM_{n+1}) (\sigma_i) &=
\bordermatrix{
 & & i & i+1 & \cr
 & I & & & \cr
 & & 0 & M_i & \cr
 & & I & I - N_i & \cr
 & & & & I
} \cdot
\Diag (TYM_{n+1} (\sigma_{i+1}), \ldots, TYM_{n+1} (\sigma_{i+1})),
\end{align*}
where
\begin{align*}
M_i &= q^2 TYM_{n+1} (\chi_n (x_{i+1})) = q^2 TYM_{n+1} ((\sigma_{i} \cdots \sigma_1)^{-1} \sigma_{i+1}^2 (\sigma_{i} \cdots \sigma_1))\\
&= q^2 \Diag (t, 1, \ldots, 1, \overset{i+2}{t}, 1, \ldots, 1)\\
\end{align*}
and
\begin{align*}
N_i &= q^2 TYM_{n+1} (\chi_n (x_{i+1}^{-1} x_i x_{i+1})) = q^2 TYM_{n+1} ((\sigma_{i} \cdots \sigma_1)^{-1} \sigma_{i+1}^{-2} \sigma_i^2 \sigma_{i+1}^2 (\sigma_{i} \cdots \sigma_1))\\
&= q^2 \Diag (t, 1, \ldots, 1, \overset{i+1}{t}, 1, \ldots, 1).
\end{align*}
It directly follows from the definitions that
\begin{align*}
&TYM_{n+1} (\sigma_{i+1}) = (1) \oplus TYM_n (\sigma_i),\\
&M_i = (q^2 t) \oplus \eta_n (x_{i+1}),\\
&N_i = (q^2 t) \oplus \eta_n (x_i) = (q^2 t) \oplus \eta_n (x_{i+1}^{-1} x_i x_{i+1}), 
\end{align*}
and thus we have
\begin{align*}
&\mathbf{LM}_q(TYM_{n+1}) (\sigma_i) \\
=&
\left((1) \oplus I_n \right)^{\oplus i-1} \oplus
\left(
\begin{array}{cc}
0 & (q^2 t) \oplus \eta_n (x_{i+1}) \\
(1) \oplus I_n & (1-q^2 t) \oplus( I_n - \eta_n (x_{i+1}^{-1} x_i x_{i+1}))
\end{array}
\right)
\oplus \left((1) \oplus I_n) \right)^{\oplus n-i-1}\\
\quad & \cdot \Diag ((1) \oplus TYM_n (\sigma_i), \ldots, (1) \oplus TYM_n (\sigma_i)).
\end{align*}
Therefore, the subspace spanned by $\{\mathbf{e}_1, \mathbf{e}_{(n+1)+1}, \mathbf{e}_{2(n+1)+1}, \ldots, \mathbf{e}_{(n-1)(n+1)+1} \}$ is invariant under $\mathbf{LM}_q(TYM_{n+1})$, and equivalent to the Burau representation $Bur_{n, q^2t}$, where $\{ \mathbf{e}_i, 1 \leq i \leq n(n+1) \}$ is the canonical basis.
Moreover, its complement is also invariant and equivalent to $\mathbf{LM}_q(\eta_n)$.
Hence the representation $\mathbf{LM}_q(TYM_{n+1})$ splits, and we deduce that:
\begin{prop}
There is a $\B_n$-module isomorphism $\mathbf{LM}_q(TYM_{n+1}) \cong Bur_{n, q^2t} \oplus \mathbf{LM}_q(\eta_n)$.
\end{prop}

From this result, we are able to rewrite Question~\ref{kernel} in terms of the representation $\eta_n \colon \F_n \rtimes_{a_n} \B_n \longrightarrow \mathrm{GL}_n(\mathbb{Z}[t^{\pm1}, q^{\pm1}]))$:

\begin{q}
Is it true that $\ker (\mathbf{LM}(\eta_{n})) = \ker (Bur_n)$
and that $\ker (\mathbf{LM}_q(\eta_{n}))$ is trivial ?
\end{q}

In addition, by considering the reduced version of the Long-Moody construction \cite[Theorem 2.11]{Long}, we obtain the reduced Long-Moody construction $\widetilde{\mathbf{LM}}_q(\eta_n) \colon \B_n \rightarrow \mathrm{GL}_{n(n-1)}(\mathbb{Z}[t^{\pm1}, q^{\pm1}]))$.
For $n=3$, this is the $6$-dimensional representation defined by the following matrices in the canonical basis induced by the representation $\eta_n$ and the Long-Moody construction and that we denote by $\{\mathbf{e}_{i},1\leq i\leq6\}$.
\begin{align*}
\sigma_1 \to
\left(
\begin{array}{cccccc}
 0 & -q^2 & 0 & 0 & 0 & 0 \\
 -q^2 t^2 & 0 & 0 & 0 & 0 & 0 \\
 0 & 0 & -q^2 & 0 & 0 & 0 \\
 0 & 1 & 0 & 0 & 1 & 0 \\
 t & 0 & 0 & t & 0 & 0 \\
 0 & 0 & 1 & 0 & 0 & 1 \\
\end{array}
\right),\ 
\sigma_2 \to
\left(
\begin{array}{cccccc}
 1 & 0 & 0 & q^2 & 0 & 0 \\
 0 & 0 & 1 & 0 & 0 & q^2 \\
 0 & t & 0 & 0 & q^2 t^2 & 0 \\
 0 & 0 & 0 & -q^2 & 0 & 0 \\
 0 & 0 & 0 & 0 & 0 & -q^2 \\
 0 & 0 & 0 & 0 & -q^2 t^2 & 0 \\
\end{array}
\right).
\end{align*}
Let $ \mathbb{Q}(t,q)$ denote the field of fractions of $\mathbb{Z}[t^{\pm1}, q^{\pm1}]$ and $\overline{\mathbb{Q}(t,q)}$ denote the algebraic closure of $ \mathbb{Q}(t,q)$. Using Mathematica, we compute that:
\begin{lemma}\label{lem:spectrum_LM(TYM)}
The spectrum of $\widetilde{\mathbf{LM}}_q(\eta_3) \otimes_{\mathbb{Z}[t^{\pm1}, q^{\pm1}]} \overline{\mathbb{Q}(t,q)}$ is $\{1,\sqrt{t},-\sqrt{t},q^{2}t,-q^{2}t,-q^{2}\}$.
\end{lemma}
Furthermore, we are also interested in the irreducibility properties of this new representation $\widetilde{\mathbf{LM}}_q(\eta_n)$. We recall that, for a group $G$ and an integral domain $R$ with field of fractions $K$, a representation $\rho:G\to \mathrm{GL}_{n}(R)$ is said to be \emph{absolutely irreducible} if $\rho\otimes_{R}\bar{K}$ is irreducible (where $\bar{K}$ denotes the algebraic closure of $K$). In particular, if such $\rho$ is absolutely irreducible, then $\rho\otimes_{R}K'$ is irreducible for $K'$ any field extension of $K$. We prove that:
\begin{thm}\label{thm:irred}
For each natural number $3\leq n\leq 6$, the $\B_n$-representation $\widetilde{\mathbf{LM}}_q(\eta_n)$ is absolutely irreducible.
\end{thm}
Preliminarily, the proof of Theorem~\ref{thm:irred} requires the following general result on the connection between irreducibility and specialisation of representations: 
\begin{lemma}\label{lem:absolute_irred_general}
Using the above notations, we consider a representation $\rho:G\to \mathrm{GL}_{n}(R)\subset\mathrm{GL}_{n}(K)$ and a ring homomorphism $\phi:R\to\mathbb{K}$ with $\mathbb{K}$ an algebraically closed field.
If the representation $\rho':G\to\mathrm{GL}_{n}(\mathbb{K})$ induced by $\phi$ is irreducible, then $\rho$ is absolutely irreducible.
\end{lemma}
\begin{proof}
Up to replacing $R$ by its localisation at $\mathrm{Ker}(\phi)$, we may assume without loss of generality that $R$ is a local ring.
As a consequence of Burnside's theorem (see for instance \cite[Theorem 2 and Section 2]{Lam}), the irreducibility of $\rho'$ implies that the image of the group algebra $\mathbb{K}[G]$ inside the matrix algebra $\mathrm{Mat}_{n}(\mathbb{K})$ is of dimension $n^2$.
Since $R$ is a local ring, as a corollary of Nakayama's lemma (see for instance \cite[Section I.3]{Lang}), we know that the $n^2$ elements of $G$ which generate $\mathrm{Mat}_{n}(\mathbb{K})$ also generate $\mathrm{Mat}_{n}(R)$ as a $R$-module. A fortiori, the image of the group algebra $\bar{K}[G]$ inside $\mathrm{Mat}_{n}(\bar{K})$ of dimension $n^2$, and thus $\rho\otimes_{R}\bar{K}$ is irreducible, which ends the proof.
\end{proof}

\begin{proof}[Proof of Theorem~\ref{thm:irred}]
It follows from Lemma \ref{lem:absolute_irred_general} that it is enough to prove that some specialized versions of $\widetilde{\mathbf{LM}}_q(\eta_n)$ are irreducible.
Namely, we assign particular integral values to the variables $t$ and $q$ by considering a certain map $\mathbb{Z}[t^{\pm1},q^{\pm1}]\to\mathbb{Z}$, and we denote the obtained representation $\widetilde{\mathbf{LM}}_{q}(\eta_{n})\otimes_{\mathbb{Z}[t^{\pm1},q^{\pm1}]}\mathbb{Z}$ by $\widetilde{\mathbf{LM}}_{q}(\eta_{n})^{\dagger}$. 
Then we run GAP codes\footnote{See \url{https://www.dropbox.com/sh/nh9d2uf27jj95al/AAA_yKl6ELeK8Gtw-FUM_bFea?dl=0}} which determine the dimension of a (non-trivial) subrepresentation of $\widetilde{\mathbf{LM}}_{q}(\eta_{n})^{\dagger}\otimes_{\mathbb{Z}}\mathbb{C}$ for $3\leq n\leq 6$ for random values of $t$ and $q$.
The obtained dimension is the one of the total space for at least one specialization: this works for the specializations $q :=3; t :=5$ and $q :=6; t :=10$ for any $3\leq n\leq 6$. Therefore the representation $\widetilde{\mathbf{LM}}_{q}(\eta_{n})^{\dagger}\otimes_{\mathbb{Z}}\mathbb{C}$ is irreducible by Burnside's theorem, which ends the proof.
\end{proof}
In particular, we directly deduce the following property for the representation over $\mathbb{Z}[t^{\pm1}, q^{\pm1}]$.
\begin{cor}
For $n\leq6$, the $\B_n$-representations $\widetilde{\mathbf{LM}}_q(\eta_n)$ with ground ring $\mathbb{Z}[t^{\pm1},q^{\pm1}]$ are indecomposable.
\end{cor}
\begin{proof}
Let us assume that $\widetilde{\mathbf{LM}}_q(\eta_n)$ is isomorphic to a direct sum $W\oplus W'$ of $\mathbb{Z}[t^{\pm1},q^{\pm1}][\B_n]$-modules.
Since $\widetilde{\mathbf{LM}}_q(\eta_n)\otimes_{\mathbb{Z}[t^{\pm1}, q^{\pm1}]} \mathbb{Q}(t,q)$ is irreducible, either $W\otimes_{\mathbb{Z}[t^{\pm1}, q^{\pm1}]} \mathbb{Q}(t,q)$ or $W'\otimes_{\mathbb{Z}[t^{\pm1}, q^{\pm1}]} \mathbb{Q}(t,q)$ is trivial.
Recall that for a $\mathbb{Z}[t^{\pm1}, q^{\pm1}]$-module $M$, the kernel of the map $M\rightarrow M\otimes_{\mathbb{Z}[t^{\pm1}, q^{\pm1}]} \mathbb{Q}(t,q)$ is the torsion submodule of $M$.
Therefore either $W$ or $W'$ is trivial since they are both free as $\mathbb{Z}[t^{\pm1}, q^{\pm1}]$-modules, which ends the proof.
\end{proof}
More generally, we make the following conjecture for the further representations for higher $n$:
\begin{conj}
The $\B_n$-representation $\widetilde{\mathbf{LM}}_q(\eta_n)$ is absolutely irreducible (and thus indecomposable) for all $n\geq7$.
\end{conj}

Finally, we highlight the newness of the representation $\widetilde{\mathbf{LM}}_q(\eta_n)$ for $n=3$ by comparing it to the well-known families of representations, namely the Burau representation, the Tong-Yang-Ma representation and the \emph{Lawrence-Krammer-Bigelow representation} $\mathbf{LKB}_{3}:\B_{3}\rightarrow \mathrm{GL}_{3}(\mathbb{Z}[t^{\pm1},q^{\pm1}])$.
We recall that this last representation has been introduced by Lawrence \cite{Lawrence1}, Bigelow \cite{bigelow2001braid} and Krammer \cite{KrammerLK} in different ways, and that Bigelow \cite{bigelow2001braid} and Krammer \cite{KrammerLK} independently proved to be faithful.
Since the dimensions of these representations are less than the one $\widetilde{\mathbf{LM}}_q(\eta_3)$, it only makes sense to wonder whether or not one of them is "contained" in $\widetilde{\mathbf{LM}}_q(\eta_3)$ or else if tensor powers of these representations are isomorphic $\widetilde{\mathbf{LM}}_q(\eta_3)$. 
However, we prove that:

\begin{prop}\label{prop:not_equiv_LKB}
The $\B_3$-representations $Bur_{3}\otimes_{\mathbb{Z}[t^{\pm1}]}\mathbb{Z}[t^{\pm1},q^{\pm1}]$, $TYM_{3}\otimes_{\mathbb{Z}[t^{\pm1}]}\mathbb{Z}[t^{\pm1},q^{\pm1}]$ and $\mathbf{LKB}_{3}$ are neither subrepresentations or quotients of $\widetilde{\mathbf{LM}}_q(\eta_3)$.

Moreover, the $\B_3$-representations $\overline{Bur}_{3,t} \otimes TYM_{3,q}$, $\overline{Bur}_{3,q} \otimes TYM_{3,t}$, $\overline{Bur}_{3,t} \otimes LKB_{3,t,q}$ and $\overline{Bur}_{3,q} \otimes LKB_{3,t,q}$ are not isomorphic to $\widetilde{\mathbf{LM}}_{q}(\eta_{3})$.
\end{prop}
\begin{proof}
Since the spaces of the representations are free $\mathbb{Z}[t^{\pm1},q^{\pm1}]$-modules, the first part of the results straightforwardly follows from the irreducibility of $\widetilde{\mathbf{LM}}_q(\eta_3) \otimes_{\mathbb{Z}[t^{\pm1}, q^{\pm1}]} \mathbb{Q}(t,q)$ by Theorem~\ref{thm:irred} together with the exactness properties of the tensor product functor $- \otimes_{\mathbb{Z}[t^{\pm1}, q^{\pm1}]} \mathbb{Q}(t,q)$.

The second part of the results follows from the comparison of the spectrum of the tensor products of the representations computed in Section~\ref{s:spectrum of representations} Lemma to the spectrum of $\widetilde{\mathbf{LM}}_{q}(\eta_{3})$ computed in Lemma~\ref{lem:spectrum_LM(TYM)}.
\end{proof}

\section*{Acknowledgement}
The authors wish to thank Takuya Sakasai who is Takano’s supervisor for introducing us to these questions. They also would like to thank Paolo Bellingeri and Tetsuya Ito for their comments and questions.
Finally, they thank the anonymous referee for their careful reading and suggestions.
The first author was supported by the Institute for Basic Science IBS-R003-D1, by a Rankin-Sneddon Research Fellowship of the University of Glasgow and by the ANR Project AlMaRe ANR-19-CE40-0001-01.

\section{Appendix A: spectrum}\label{s:spectrum of representations}
We consider the Tong-Yang-Ma representation (with $\mathbb{Z}[t^{\pm1},q^{\pm1}]$ as ground ring) $TYM_{3,t} \colon B_3 \longrightarrow GL_3(\mathbb{Z}[t^{\pm1},q^{\pm1}])$ is defined by
$$\sigma_1 \longmapsto
\left(
\begin{array}{ccc}
0 & 1 & 0 \\
t & 0 & 0 \\
0 & 0 & 1 \\
\end{array}
\right),\ 
\sigma_2 \longmapsto
\left(
\begin{array}{ccc}
1 & 0 & 0 \\
0 & 0 & t \\
0 & 1 & 0 \\
\end{array}
\right),$$
that the \emph{reduced} Burau representation (with $\mathbb{Z}[t^{\pm1},q^{\pm1}]$ as ground ring) $\overline{Bur}_{3,t} \colon B_3 \longrightarrow GL_2(\mathbb{Z}[t^{\pm1},q^{\pm1}])$ is defined by
$$\sigma_1 \longmapsto
\left(
\begin{array}{cc}
-t & 1 \\
0 & 1 \\
\end{array}
\right),\ 
\sigma_2 \longmapsto
\left(
\begin{array}{cc}
1 & 0 \\
t & -t \\
\end{array}
\right),$$
and that the Lawrence-Krammer-Bigelow representation $LKB_{3,t,q} \colon B_3 \longrightarrow GL_3(\mathbb{Z}[t^{\pm1},q^{\pm1}])$ (see Krammer \cite[Section 3]{KrammerLK} or Paoluzzi and Paris \cite{PaoluzziParis}) is defined by
$$\sigma_1 \longmapsto
\left(
\begin{array}{ccc}
 -q^2 t & q(1-q)t & 0 \\
 0 & 0 & q \\
 0 & 1 & 1-q \\
\end{array}
\right),\ 
\sigma_2 \longmapsto
\left(
\begin{array}{ccc}
 1-q & 1 & 0 \\
 q & 0 & 0 \\
 0 & q(1-q)t & -q^2 t \\
\end{array}
\right).$$

Let $A$ and $B$ be square matrices of size $n$ and $m$, and  $\lambda_1, \ldots, \lambda_n$ and $\mu_1, \ldots, \mu_m$ the eigenvalues of $A$ and $B$, respectively.
Then the eigenvalues of $A \otimes B$ are $\lambda_i \mu_j$ $(1 \leq i \leq n, 1 \leq j \leq m)$.
The spectrums of $TYM_{3,t}, \overline{Bur}_{3,t}$ and $LKB_{3,t,q}$ over $\overline{\mathbb{Q}(t,q)}$ are $\{ 1, \sqrt{t}, -\sqrt{t} \}, \{ 1, -t \}$ and $\{ 1, -q, -q^2t \}$, respectively.
Hence, the spectrums of 6-dimensional representations $\rho \colon B_3 \longrightarrow GL_6(\Z[t^{\pm1}, q^{\pm1}])$ given by three representations are
$$
\begin{cases}
\{1,-\sqrt{q},\sqrt{q},-t,-\sqrt{q} t,\sqrt{q} t\} & $if$\ \rho = \overline{Bur}_{3,t} \otimes TYM_{3,q},\\
\{1,-q,-\sqrt{t},\sqrt{t},-q \sqrt{t},q \sqrt{t}\} & $if$\ \rho = \overline{Bur}_{3,q} \otimes TYM_{3,t},\\
\{1,-q,-t,q t,-q^2 t,q^2 t^2\} & $if$\ \rho = \overline{Bur}_{3,t} \otimes LKB_{3,t,q},\\
\{1,-q,-q,q^2,-q^2 t,q^3 t\} & $if$\ \rho = \overline{Bur}_{3,q} \otimes LKB_{3,t,q}.
\end{cases}
$$

\bibliography{draft}
\bibliographystyle{plain}

\end{document}